\documentclass[11pt]{article}

\usepackage{a4wide}

\usepackage{parskip}
\usepackage[english]{babel}

\usepackage{graphicx}

\usepackage{amssymb}

\usepackage{amsmath}

\usepackage{xcolor}
\usepackage{url} % Puse este paquete para escribir las direcciones web de los paquetes de GAP

\begin{document}

%%%%%%%%%%%%%%%%%%%%%%%%%%%%%%%%%%%%%%%%%%%%
%DEFINICIONES NUESTRAS
%%%%%%%%%%%%%%%%%%%%%%%%%%%%%%%%%%%%%%%%%%%%
\newcommand{\Nat}{{\mathbb N}}
\newcommand{\Real}{{\mathbb R}}
\newcommand{\Z}{{\mathbb Z}}
\newcommand{\noS}{\overline{S}}
\newcommand{\f}{\mathfrak{f}}
\newcommand{\p}{\mathrm{p}}
\newcommand{\wa}{W_a}
\newcommand{\wb}{W_b}
\newcommand{\mcd}{\mathrm{mcd}}
\newcommand{\dg}{\mathrm{d}_G}
\newcommand{\dd}{\mathrm{d}}
\newcommand{\EE}{\mathcal{E}}
\newcommand{\F}{\mathcal{F}}
\newcommand{\T}{\mathcal{T}}
\newcommand{\HH}{\mathcal{H}}
\newcommand{\LL}{\mathcal{L}}
\newcommand{\ele}{\mathrm{L}}
\newcommand{\Ap}{\mathrm{Ap}}
\newcommand{\Su}{\mathbb{S}}
\def\vecu{\mbox{\boldmath$u$}}
\def\vecv{\mbox{\boldmath$v$}}
\def\vece{\mbox{\boldmath$e$}}
\def\x{\mathrm{x}}
\def\y{\mathrm{y}}
\def\z{\mathrm{z}}
\newcommand{\sg}[1]{\langle{#1}\rangle}
\newcommand{\sq}[1]{[\![{#1}]\!]}
\newcommand{\mmod}[1]{\!\!\pmod{#1}}
\newcommand{\lcm}{\mathrm{lcm}}
\newcommand{\w}{\mathrm{W}}
\newcommand{\PP}{\mathcal{P}}
\newcommand{\bb}[1]{\partial{#1}}
%%%%%%%%%%%%%%%%%%%%%%%%%%%%%%%%%%%%%%%%%%%%%%%
%%%%%%%%%%%%%%%%%%%%%%%%%%%%%%%%%%%%%%%%%%%%%%%
\newtheorem{lem}{Lemma}
\newtheorem{pro}{Proposition}
\newtheorem{defi}{Definition}
\newtheorem{teo}{Theorem}
\newtheorem{cor}{Corollary}
\newtheorem{exa}{Example}
\newtheorem{rem}{Remark}

\def\peb[#1]{{\left\lfloor #1\right\rfloor}}
\newcommand{\bs}{\boldsymbol}

\def\cardinal{\sharp}
%%%%%%%%%%%%%%%%%%%%%%%%%%%%%%%%%%%%%%%%%%%%%%%%
\newcommand{\algP}{\textsf{P}}
\newcommand{\algL}{\textsf{L}}
\newcommand{\algBCS}{\textsf{BCS}}
\newcommand{\algAL}{\textsf{AL}}

\title{Computing denumerants in\\
numerical $3$--semigroups\thanks{The second author is supported by the projects MTM2014-55367-P, FQM-343 and FEDER funds. The first author is supported by the project
MTM2014-60127-P.}}
      
\author{F. Aguil\'o--Gost \and D. Llena}
\date{January 16, 2017}

\maketitle

\begin{abstract}
As far as we know, usual computer algebra packages can not compute denumerants for almost medium (about a hundred digits) or almost medium--large (about a thousand digits) input data in a reasonably time cost on an ordinary computer. Implemented algorithms can manage numerical $n$--semigroups for small input data.

Here we are interested in denumerants of numerical $3$--semigroups which have almost medium input data. A new algorithm for computing denumerants is given for this task. It can manage almost medium input data in the worst case and medium--large or even large input data in some cases. 
\end{abstract}

\vspace*{3mm}


\section{Introduction}

Let $\Nat$ be the set of non negative integers. We denote the equivalence class of $k$ modulo $m$ as $[k]_m$. Given $n_1,\ldots,s_k\in\Nat$, $1<n_1<\cdots<n_k$ and $\gcd(n_1,\ldots,n_k)=1$, the \textit{numerical $k$--semigroup} $T$ generated by $G=\{n_1,\ldots,n_k\}$ is defined by
\[
T=\sg{n_1,\ldots,n_k}=\{x_1n_1+\cdots+x_kn_k:~x_1,\ldots,x_k\in\Nat\}.
\]
The {\em generating set} $G$ has not necessarily be minimal. The cardinality of a minimal generating set is the {\em embedding dimension}, $e(T)$, of the semigroup. Given an element $m\in T\setminus\{0\}$, the \textit{Ap\'ery} set of $T$ with respect to $m$ is the set $\Ap(m,T)=\{s\in T:~s-m\notin T\}$. It is well known the equivalence $s\in\Ap(m,T)\Leftrightarrow s=\min([s]_m\cap T)$ and so, $\Ap(m,T)=\{s_0,\ldots,s_{m-1}\}$ with $s_i\equiv i\mmod{m}$.

Given $s\in T$, a vector $(x_1,\ldots,x_n)\in\Nat^k$ such that $x_1n_1+\ldots+x_kn_k=s$ is called a \textit{factorization} of $s$ in $T$. Let us denote the set of factorizations of $s$ in $T$ by
\[
\F(s,T)=\{(x_1,\ldots,x_k)\in\Nat^k:~x_1s_1+\cdots+x_ks_k=s\}.
\]
The \textit{denumerant of $s$ in $T$} is defined as the cardinality of the set $\F(s,T)$, denoted by $\dd(s,T)=|\F(s,T)|$. The {\em Frobenius number} of $T$ is defined by $\f(T)=\max(\Nat\setminus T)$. Detailed results on numerical semigroups can be found in the book of J. C. Rosales and P. A. Garc\'{\i}a--S\'anchez \cite{RoGa-NS:2009}. It is also interesting the book of Ram\'{\i}rez Alfons\'{\i}n \cite{Ra:2005} where it can be found a complet source of results related to Frobenius number.

Sylvester \cite{Sy:1882} in 1882 gave the generating function $\phi(z)$ of $\dd(m,\sg{n_1,\ldots,n_k})$
\[
\phi(z)=\frac1{(1-z^{n_1})(1-z^{n_2})\cdots(1-z^{n_k})}.
\]
Schur \cite{Sch:1926} in 1926 studied the asymptotic behaviour of the denumerant,
\[
\limsup_{m\to\infty}\frac{\dd(m,\sg{n_1,\ldots,n_k})}{\frac{m^{k-1}}{n_1\cdots n_k(k-1)!}}=1.
\]
Sylvester \cite{Sy:1857} in 1857 and Cayley \cite{Ca:1860} in 1860 gave the expression $\dd(m,\sg{n_1,\ldots,n_k})=P_k(m)+Q_k(m)$ where $P_k(m)$ is a polynomial of degree $k-1$ and $Q_k(m)$ is a periodic function in the variable $m$. Beck, Gessel and Komatsu \cite{BeGeKo:2001} in 2001 found an expression for $P_k(m)$ that depends upon Bernoulli numbers.

Popoviciu \cite{Po:1953} in 1953 found an efficient semi--closed expression\footnote{This expression only requires $O(\log\max\{p,q\})$ arithmetic operations to be applied.} for $n=2$ of $\dd(m,\sg{p,q})$
\[
\dd(m,\sg{p,q})=\frac{m+pf(m)+qg(m)}{\lcm(p,q)}-1,
\]
where $f(m)\equiv-mp^{-1}\mmod{q}$ with $1\leq f(m)\leq q$ and $g(m)\equiv-mq^{-1}\mmod{p}$ with $1\leq g(m)\leq p$.  Ehrhart \cite{EhII:1967} in 1967 and Sert\"oz and \"Ozl\"uk in 1991 gave recursive denumerant formulae for $2\leq k\leq4$. You can find an exhaustive set of results on denumerants in the book of J. Ram\'{\i}rez Alfons\'{\i}n \cite{Ra:2005}.

No similar efficient semi--closed expressions are known for $k\geq3$, however there are some known numerical algorithms to find the set of factorizations $\F(m,T)$ in the general case. Unfortunately, as far as we know, usual computer algebra systems have implemented no command for denumerant. Thus, the calculation of denumerant turns to be a time consuming task. Taking for instance, $n_1=7^k$, $n_2=11^k$, $n_2=\f(\sg{7^k,11^k})$, $P_k=n_1n_2n_3$, $S_k=n_1+n_2+n_3$ and $m_k=P_k-S_k-k$, we obtain the figures of Table~\ref{tab:cascomp} for $\dd(m_k,T_k)$ and $T_k=\sg{n_1,n_2,n_3}$. The reason why we choose $m_k$ is clear by Theorem~\ref{teo:SeOz}.

\begin{table}[h]
\centering
%\hspace*{-0.7cm}
%\small
\footnotesize
\begin{tabular}{|r|r|r|r|r|r|}
\hline
%\multicolumn{7}{|c|}{$T_k=\sg{7^k,11^k,\mathfrak{f}(7^k,11^k)},\quad m_k=P_k-S_k-k$} \\
%\hline
$k$&$m_k$&$\dd(m_k,T_k)$&{\sf Mathematica~8}&{\sf Sage~7.3}&{\sf GAP~1.5.1}\\
\hline\hline
$1$&$4465$&$2232$&0.011311&0.019617&0.009835\\
$2$&$34139180$&$17069589$&177.318173&535.270590&6.100101\\
\hline
\end{tabular}
\caption{Time in seconds using an \textsf{i5}@1.3Ghz processor}
\label{tab:cascomp}
\end{table}

Table~\ref{tab:cascomp} shows how popular CAS programs\footnote{The commands for computing the denumerant $\dd(m,\sg{a,b,c})$ are \texttt{Length[FrobeniusSolve[\{a,b,c\},m]]} for {\sf Mathematica~8}, \texttt{WeightedIntegerVectors(m,[a,b,c]).cardinality()} for {\sf Sage~7.3} and \texttt{NrRes\-tricted\-Parti\-tions(m,[a,b,c])} for {\sf GAP~1.5.1}} can not manage almost medium (about half a hundred digits). Clearly the {\sf Gap} package takes advantage for these input instances. From now on, we focus our attention to denumerants of numerical $3$--semigroups and the notation $n_1=a$, $n_2=b$, $n_3=c$ and $T=\sg{a,b,c}$ will be used here.

Popoviciu \cite[page 27]{Po:1953} gave an $O(c\log c)$ algorithm, in the worst case, for computing $\dd(m,T)$ when $\{a,b,c\}$ are pairwise coprime numbers (pcn). Lison\v{e}k \cite[page 230]{Li:1995} in 1995 gave an $O(ab\log b)$ algorithm, in the worst case for pcn (this time cost  can be reduced to $O(ab)$ provided that a number of $\max\{O(a^2b^2),\allowbreak O(abc)\}$ precomputed values, related to $T$, can be stored in the computer memory for later usage). Brown, Chou and Shiue in 2003 \cite[page 199]{BCS:2003} gave an $O(ab\log c)$ algorithm, in the worst case. This last work also contains interesting results on denumerants that can be taken into account for numerical calculations. We refer to these algorithms as \algP, \algL\  and \algBCS, respectively. Notice that the speed of Algorithm~\algP\ versus Algorithm~\algL\ depends on the ratio $\frac{c\log c}{ab\log b}$.

Algorithms \algP, \algL\ and \algBCS\ calculate the denumerants of Table~\ref{tab:cascomp} significantly faster. A non-compiled {\sf Sage~7.3} implementations of them give the figures in Table~\ref{tab:algs} (using the same processor of Table~\ref{tab:cascomp}).

\begin{table}[h]
\centering
%\footnotesize
\small
\begin{tabular}{|r|r|r|r|r|r|}
\hline
$k$&\multicolumn{1}{|c}{$m_k$}&\multicolumn{1}{|c}{$\dd(m_k,T_k)$}&\multicolumn{1}{|c}{\algP} &\multicolumn{1}{|c}{\algL} &\multicolumn{1}{|c|}{\algBCS} \\
\hline\hline
1&4465&2232&0.003994&0.004331&0.011661\\
2&34139180&17069589&0.083913&0.138889&0.920211\\
3&207657687311&103828843654&5.063864&9.495915&63.647251\\
\hline
\end{tabular}
\caption{Time in seconds obtained by \algP, \algL\ and \algBCS}
\label{tab:algs}
\end{table}

Nonetheless, these algorithms do not reach the necessary efficiency for managing almost medium input. The goal of this work is to provide a reasonably efficient new algorithm which allows such kind of inputs when working on ordinary computers.

Our algorithm has a theoretical time cost of $O(b+\log c)$, in the worst case. However, numerical evidences suggest that, in some cases, it can have a smaller cost\footnote{As an instance, the same data of Table~\ref{tab:algs} for $k=10^3$ is calculated in $0.009836$ seconds and for $k=10^5$ in $2.110464$ seconds.}. This algorithm is based on a semi-closed denumerant expression given in \cite{AG:2010} which is included here in Theorem~\ref{teo:BS}.

The summary of the paper is the following: Section~$2$ contains the basic known tools, mainly Theorem~\ref{teo:BS} and expression \eqref{eq:BS}. Section~$3$ developes expression \eqref{eq:BS} to be used for numerical purposes. In this developing it is apparent that the main computation depends on the so called $S^\pm$ {\em discrete sums}. Some tools to calculate $S^\pm$ sums are developed in Section~$4$, mainly the so called {\em hS-type sets}. Section~$5$ contains the main algorithm and Section~$6$ analyzes the time cost, in the worst case. Finally, in Section~$7$, several instances of time tests are given.

\section{Some definitions and known results}

In this section we give the main known results that allow us to reach our goal. The usual notation for semigroups will be $T=\sg{a,b,c}$ with $1\leq a<b<c$ and $\gcd(a,b,c)=1$. Also the product $P=abc$ and sum $S=a+b+c$ of the generators are used.

Although algorithms {\sf P} and {\sf L} act over pairwise coprime generators, this condition can be removed by the following result due to Brown, Shou and Shiue \cite{BCS:2003}. Here the integer $u'_v(t)$ is defined to be the unique integer value $1\leq u'_v(t)\leq v$ such that $uu'_v(t)\equiv-t\mmod{v}$ with $u,v\geq1$ and $\gcd(u,v)=1$.

\begin{lem}[Brown, Chou and Shiue 2003 {\rm\cite[Lemma~4.5]{BCS:2003}}]\label{lem:coprim}
Consider the semigroup $T=\sg{a,b,c}$ with $\gcd(a,b,c)=1$. Set $g_a=\gcd(b,c)$, $g_b=\gcd(a,c)$ and $g_c=\gcd(a,b)$. For any integer $n>0$, the integer value $n'=n-(g_a-a'_{g_a}(n))a-(g_b-b'_{g_b}(n))b-(g_c-c'_{g_c}(n))c$ is multiple of $g_ag_bg_c$ and the denumerant's identity $\dd(n,T)=\dd(\frac{n'}{g_ag_bg_c},T')$ holds with $T'=\sg{\frac{a}{g_bg_c},\frac b{g_ag_c},\frac c{g_ag_b}}$. Here it is understood that $\dd(0,T')=1$ and $\dd(\frac{n'}{g_ag_bg_c},T')=0$ whenever $n'<0$.
\end{lem}

%COMENTARI: P--periodicitat
By the following theorem, due to Ehrhart in 1967, we only need to compute denumerants in the range of values $m\in\{0,1,\ldots,P-1\}$.

\begin{teo}[Ehrhart~1967 {\rm\cite[Theorem~10.5]{EhII:1967}}]\label{teo:EhII}
Consider $T=\sg{a,b,c}$ with  $a$, $b$ and $c$ pcn. Set $P=abc$, $S=a+b+c$ and $m=qP+r$ with $0\leq r<P$. Then, 
\[
\dd(m,T)=\dd(r,T)+\frac{q(m+r+S)}2.
\]
In particular,
\[
\dd(P,T)=\frac{P+S}2+1.
\]
\end{teo}

The range $\{0,\ldots,P-1\}$ can be reduced to $\{0,\ldots,P-S\}$ by the following theorem due to Sert\"oz and \"Ozluk in 1991.

\begin{teo}[{Sert\"oz and \"Ozl\"uk~1991 \rm\cite[page~4]{SeOz:1991}}]\label{teo:SeOz}
Consider $T=\sg{a,b,c}$ with $a$, $b$ and $c$ pcn. Set $P=abc$ and $S=a+b+c$. Then, for $1\leq x\leq S-1$ we have
\[
\dd(P-x,T)=\frac{P+S}2-x.
\]
In particular,
\[
\dd(P-S+1)=\frac{P-S}2+1.
\]
\end{teo}

\begin{rem}\label{rem:tcomp0}\sl
The time cost, in the worst case, of the algorithms {\sf P}, {\sf L} and {\sf BCS} for computing the denumerant $\dd(m,\sg{a,b,c})$ have been given for the largest value of $m$ (by theorems \ref{teo:EhII} and \ref{teo:SeOz}), that is $m\approx P=abc$.
\end{rem}

%COMENTARI: L--formes associades a un semigrup.
We use the concept of {\em L--shape} as a main tool for the new algorithm. Thus, we include here some known results for this geometrical discrete structure. Denote the {\em interval} $[s,t)=\{x\in\Real:~s\leq x<t\}$, the {\em unitary square} $\sq{m,n}=[m,m+1)\times[n,n+1)\in\Real^2$ and the {\em discrete backwards cone} $\Delta(u,v)=\{\sq{m,n}:~(m,n)\in\Nat^2,0\leq m\leq u, 0\leq n\leq v\}$ for each $u,v\in\Nat$. We also denote the equivalence class of $u$ modulo $v$ by $[u]_v$.

Consider each unitary square $\sq{m,n}$, for $(m,n)\in\Nat^2$, labelled by the equivalence class $[ma+nb]_c$. Define the minimum values
\begin{equation}
M_n=\min\{sa+tb:~(s,t)\in\Nat^2, [sa+tb]=[n]_c\}.\label{eq:Mn}
\end{equation}

\begin{defi}[Minimum distance diagram]\label{def:mdd}
Consider a numerical $3$--semigroup $T=\langle a,\linebreak[0] b,\linebreak[0] c\rangle$. A minimum distance diagram (MDD), $\HH$, related to $T$ is a set of $c$ unitary squares that fulfils the following properties
\begin{itemize}
\item[(a)] for each $n\in\{0,\ldots,c-1\}$, there is some unitary square $\sq{s,t}\in\HH$ such that $[sa+tb]_c=[n]_c$,
\item[(b)] $\Delta(s,t)\subseteq\HH$ for each $\sq{s,t}\in\HH$,
\item[(c)] if $\sq{s,t}\in\HH$, then $sa+tb=M_n$ with $[sa+tb]_c=[n]_c$ and $M_n$ defined by \eqref{eq:Mn}.
\end{itemize}
\end{defi}

\begin{figure}[h]
\centering
\includegraphics[width=0.5\linewidth]{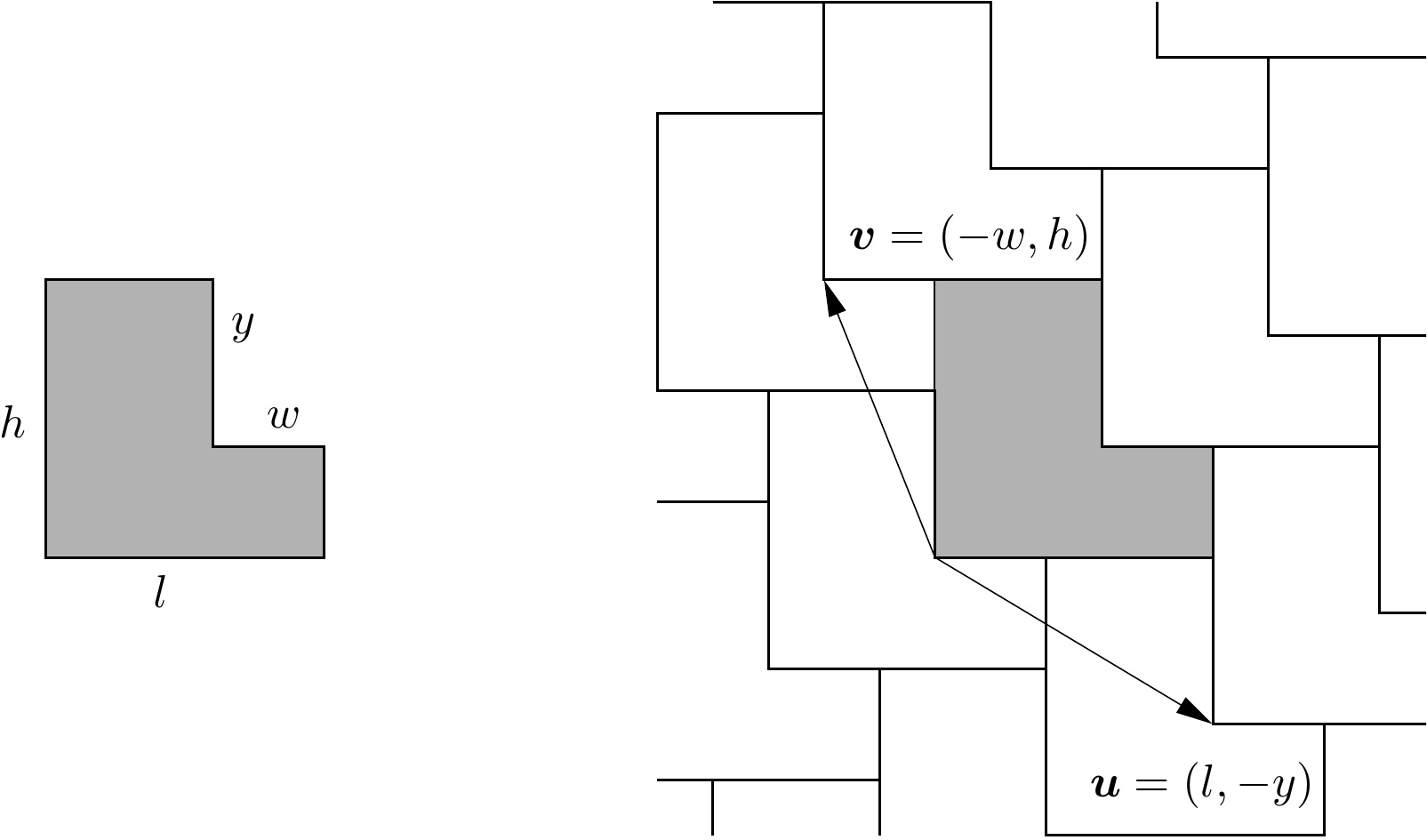}
\caption{Generic L-shape and its related tessellation}
\label{fig:eletes}
\end{figure}

Minimum distance diagrams related to numerical $3$--semigroups are known to be L--shapes or rectangles (that will be considered as degenerated L-shapes). For this reason we also refer to MDD as L-shapes and they are denoted by the lengths of their sides $\ele(l,h,w,y)$, see Figure~\ref{fig:eletes}, with $0\leq w<l$, $0\leq y<h$ and $lh-wy=c$. An L--shape tessellates the plane by translation through the vectors $\vecu=(l,-y)$ and $\vecv=(-w,h)$. The following result characterizes the L-shapes related to $T=\sg{a,b,c}$. From now on we assume $0<a<b<c$ and $\gcd(a,b,c)=1$.

\begin{teo}[A. and Mariju\'an 2014 \rm\cite{AM:2014}]\label{teo:LMDD}
Consider the numerical $3$--semigroup $T=\sg{a,b,c}$. An L-shape $\HH=\ele(l,h,w,y)$ is related to $T$ if and only if
\begin{itemize}
\item[(a)] $lh-wy=c$ and $\gcd(l,h,w,y)=1$,
\item[(b)] $la-yb\equiv0\mmod{c}$ and $hb-wa\equiv0\mmod{c}$,
\item[(c)] $la-yb\geq0$, $hb-wa\geq0$ and both expressions can't vanish at the same time.
\end{itemize}
\end{teo}

Each numerical $3$--semigroup has two related L-shapes at most (either one if $(la-yb)(hb-wa)>0$ or two whenever $(la-yb)(hb-wa)=0$, see \cite[theorems 2 and 3]{AM:2014}). L-shapes contain main information of the related semigroup. For instance, if a semigroup $T=\sg{a,b,c}$ has related the L-shape $\HH$, we have $\Ap(c,T)=\{ia+jb:~\sq{i,j}\in\HH\}$.

A classification of $3$--semigroups was given in terms of its related L--shapes in \cite{AM:2014}. The tessellation of the plane associated with each L--shape was used to derive the semi-closed expression \eqref{eq:BS} for the denumerant in \cite{AG:2010}.

%COMENTARI: Formula denumerant com a suma.
Given $T=\sg{a,b,c}$ and a related L--shape $\HH=\ele(l,h,w,y)$, let us denote $\delta=(la-yb)/c$ and $\theta=(hb-wa)/c$. From the definition of $\HH$ and Theorem~\ref{teo:LMDD}, it follows that
\begin{itemize}
\item $a=h\delta+y\theta$ and $b=w\delta+l\theta$,
\item $\delta,\theta\in\Nat$ and $\delta+\theta>0$,
\item $\delta=0\Rightarrow y>0$ and $\theta=\frac bl=\frac ay>0$,
\item $\theta=0\Rightarrow w>0$ and $\delta=\frac ah=\frac bw>0$,
\item $w=0\Rightarrow\theta=\frac bl>0$,
\item $y=0\Rightarrow\delta=\frac ah>0>0$.
\end{itemize}
All these properties will be used along this work.

Given $m\in T$, it is called the {\em basic factorization of $m$ with respect to $\HH$}, $(x_0,y_0,z_0)\in\F(m,T)$, the unique factorization such that $\sq{x_0,y_0}\in\HH$. This factorization can be computed in time cost $O(\log c)$ \cite{AB:2008}.

\begin{teo}[A. and P.A. Garc\'{\i}a S\'anchez 2010 \rm\cite{AG:2010}]\label{teo:BS}
Given $T=\sg{a,b,c}$ and a related L--shape $\HH=\ele(l,h,w,y)$, assume $m\in T$. Define $A_m=\left\lfloor\frac{z_0}{\delta+\theta}\right\rfloor$, where $(x_0,y_0,z_0)$ is the basic factorization of $m$ wrt $\HH$. For each $0\leq k\leq A_m$, set
\begin{equation}
S_k=\left\{\begin{array}{ll}
\left\lfloor\frac{y_0+k(h-y)}y\right\rfloor& \textrm{if }\delta=0,\\[10pt]
\left\lfloor\frac{z_0-k(\delta+\theta)}\delta\right\rfloor& \textrm{if }y=0,\\[10pt]
\min\{\left\lfloor\frac{y_0+k(h-y)}y\right\rfloor,\left\lfloor\frac{z_0-k(\delta+\theta)}\delta\right\rfloor\}&\textrm{if }\delta y\neq0,
\end{array}\right.\label{eq:Sk}
\end{equation}
and
\begin{equation}
T_k=\left\{\begin{array}{ll}
\left\lfloor\frac{x_0+k(l-w)}w\right\rfloor& \textrm{if }\theta=0,\\[10pt]
\left\lfloor\frac{z_0-k(\delta+\theta)}\theta\right\rfloor& \textrm{if }w=0,\\[10pt]
\min\{\left\lfloor\frac{x_0+k(l-w)}w\right\rfloor,\left\lfloor\frac{z_0-k(\delta+\theta)}\theta\right\rfloor\}&\textrm{if }\theta w\neq0.
\end{array}\right.\label{eq:Tk}
\end{equation}
Then, the denumerant of $m$ in $T$ is
\begin{equation}
\dd(m,T)=1+A_m+\sum_{k=0}^{A_m}(S_k+T_k).\label{eq:BS}
\end{equation}
\end{teo}
The sum appearing in this theorem is known as {\em the basic sum} of the denumerant with respect to the L--shape $\HH$. The direct computation of this sum does not give an efficient algorithm for calculating the denumerant. However, as it will be seen later, a detailed analysis of this expression does it.

\begin{figure}[h]
\centering
\includegraphics[width=0.6\linewidth]{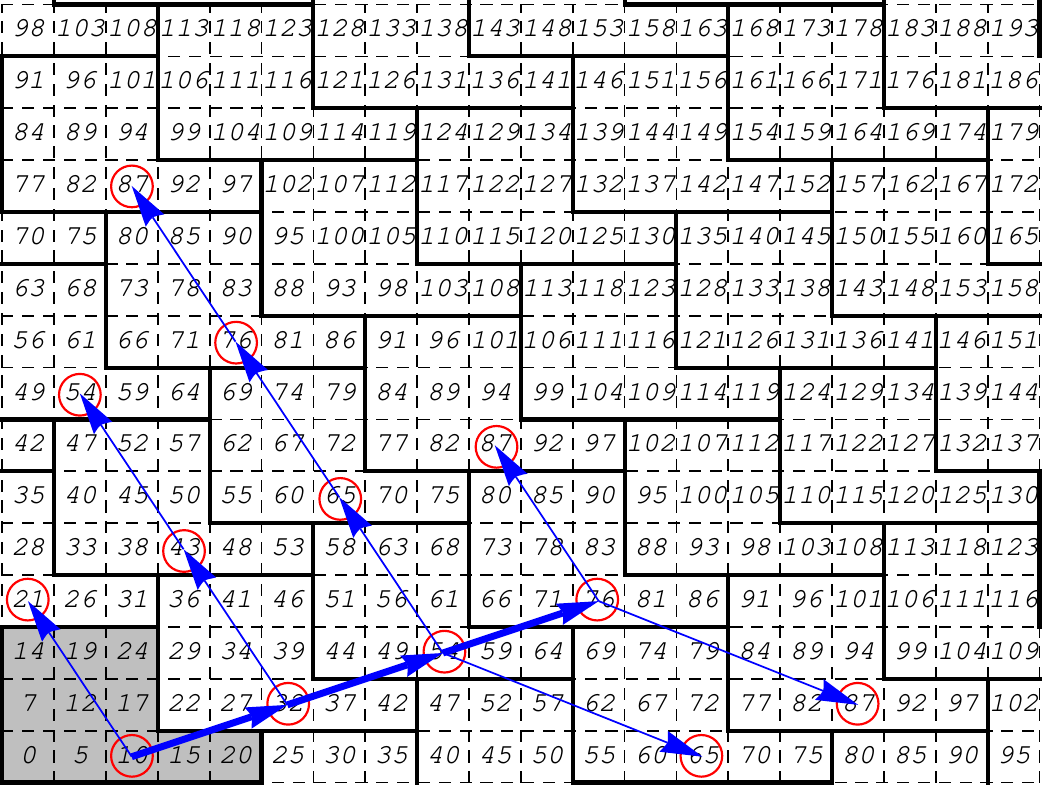}
\caption{Tree-like structure of $\pi(\F(87,\sg{5,7,11}))$}
\label{fig:tree}
\end{figure}

\begin{exa}\em\small
Take $T=\sg{5,7,11}$ and $m=87$. A related L--shape is $\HH=\ele(5,3,2,2)$, with $\delta=\theta=1$. The basic factorization of $87$ in $T$ is $(x_0,y_0,z_0)=(2,0,7)$. Thus, we have $A_m=3$. Then, it follows that $S_0+T_0=0+1=1$, $S_1+T_1=0+2=2$, $S_2+T_2=1+3=5$, $S_3+T_3=1+1=2$ and so
\[
\dd(87,T)=1+3+(1+2+5+2)=13.
\]
A geometric representation of the plane projection of the set $\F(87,T)$, $\pi(\F(87,T))$, is depicted in Figure~\ref{fig:tree}. It has a tree-like structure, given by the vectors $\vecu$, $\vecv$ and $\vecu+\vecv$ (and so, it follows the tessellation of the plane by $\HH$). Each unitary square $\sq{s,t}$ is labelled with the value $5s+7t$ (notice that values corresponding to unitary squares in the gray L--shape form the Ap\'ery set $\Ap(87,T)=\{0,12,24,14,15,5,17,7,19,20,10\}$). The unitary squares corresponding to the first two coordinates of each factorization are circled. From the coordinates of a circled unitary square $\sq{s,t}$ follows the related factorization $(s,t,\frac{87-3s-7t}{11})$. The set of factorizations is
\begin{align*}
\F(87,T)=&\{(2,0,7),(0,3,6),(5,1,5),(3,4,4),(1,7,3),(8,2,3),(13,0,2),\\
&(6,5,2),(4,8,1),(2,11,0),(11,3,1),(16,1,0),(9,6,0)\}.
\end{align*}
\end{exa}

\section{Developing the basic sum}
As it has been commented before, the basic sum \eqref{eq:BS} does not provide a direct efficient algorithm for calculating denumerants. Thus, a detailed analysis is needed. We consider three main cases: case (i) $\delta=0$, case (ii) $\theta=0$ and case (iii) $\delta\theta>0$.

The analysis of these cases reveals that the basic sum depends on several sums of the same kind. These sums will be referred to as $S^\pm$ sums and will be studied in the next section. These sums have the form $S^\pm(s,t,q,N)=\sum_{k=0}^N\left\lfloor\frac{s\pm kt}{q}\right\rfloor$ with $0\leq s,t<q$.

In this section we assume that $\HH=\ele(l,h,w,y)$ is an L--shape related to the numerical $3$--semigroup $T=\sg{a,b,c}$. We also assume that $m\in T$ and $(x_0,y_0,z_0)$ is the basic factorization of $m$ with respect to $T$.

\subsection{Case (i) $\delta=0$}

This case leads to the following expressions of the denumerant.

\begin{teo}\label{teo:casi}
Let us assume $\delta=0$. Set $A_m=\left\lfloor\frac{z_0}{\theta}\right\rfloor$. Then,
\begin{itemize}
\item[(i.1)] if $w=0$, then
\begin{equation}
\dd(m,T)=(1+A_m)(1+A_m+\overline{y_0})+(\overline{h}-2)\frac{A_m(1+A_m)}2+\sum_{k=0}^{A_m}\left\lfloor\frac{\widehat{y_0}+k\widehat{h}}{y}\right\rfloor,\label{eq:casi1}
\end{equation}
where $y_0=\overline{y_0}y+\widehat{y_0}$ with $0\leq\widehat{y_0}<y$ and $h=\overline{h}y+\widehat{h}$ with $0\leq\widehat{h}<y$.
\item[(i.2)] if $w>0$, set $k_0=\left\lceil\frac{z_0w-x_0\theta}{b}\right\rceil$. Then,
\begin{itemize}
\item[(i.2.1)] if $k_0=0$, then $\dd(m,T)$ has the same expression as in \eqref{eq:casi1}.
\item[(i.2.2)] if $1\leq k_0\leq A_m$,
\begin{align}
\dd(m,T)&=(1+A_m)(1+A_m+\overline{y_0})+k_0(\overline{x_0}-A_m)+(\overline{h}-2)\frac{A_m(1+A_m)}2\nonumber\\
&\quad+\overline{l}\frac{(k_0-1)k_0}2+\sum_{k=0}^{A_m}\left\lfloor\frac{\widehat{y_0}+k\widehat{h}}y\right\rfloor+\sum_{k=0}^{k_0-1}\left\lfloor\frac{\widehat{x_0}+k\widehat{l}}w\right\rfloor,\label{eq:casi22}
\end{align}
where $\overline{y_0},\widehat{y_0},\overline{h},\widehat{h}$ are defined as in the previous case, $x_0=\overline{x_0}w+\widehat{x_0}$ with $0\leq\widehat{x_0}<w$ and $l=\overline{l}w+\widehat{l}$ with $0\leq\widehat{l}<w$.
\item[(i.2.3)] if $k_0>A_m$,
\begin{align}
\dd(m,T)&=(1+A_m)(1+\overline{x_0}+\overline{y_0})+(\overline{l}+\overline{h}-2)\frac{A_m(1+A_m)}2\nonumber\\
&\quad+\sum_{k=0}^{A_m}\left\lfloor\frac{\widehat{x_0}+k\widehat{l}}w\right\rfloor+\sum_{k=0}^{A_m}\left\lfloor\frac{\widehat{y_0}+k\widehat{h}}y\right\rfloor,\label{eq:casi23}
\end{align}
where $\overline{y_0},\widehat{y_0},\overline{h},\widehat{h},\overline{x_0},\widehat{x_0},\overline{l}$ and $\widehat{l}$ are defined as in the previous case.
\end{itemize}
\end{itemize}
\end{teo}

\begin{rem}\sl\label{rem:k0}
Notice that $k_0\geq0$. Indeed, let us see $-1<\frac{z_0w-x_0\theta}b$ (recall that we have $w>0$). From $b=l\theta$ (recall that $\delta=0$), we have (recalling $b=l\theta$)
\[
-1<\frac{z_0w-x_0\theta}b\Leftrightarrow 0<z_0w+(l-x_0)\theta.
\]
Now, as $\sq{x_0,y_0}\in\HH$, it follows that $0\leq x_0<l$ and so the inequality $0<z_0w+(l-x_0)\theta$ holds.
\end{rem}

\noindent{\bf Proof of Theorem~\ref{teo:casi}}: If $\delta=0$, then we have $y\neq0$, $\theta>0$. From Theorem~\ref{teo:BS}, we have $A_m=\left\lfloor\frac{z_0}{\theta}\right\rfloor$ and $S_k=\left\lfloor\frac{y_0+k(h-y)}{y}\right\rfloor=\left\lfloor\frac{y_0+kh}{y}\right\rfloor-k$ for all $0\leq k\leq A_m$. Now two subcases appear, (i.1) $w=0$ and (i.2) $w>0$.
\begin{itemize}
\item[(i.1)] Assume $w=0$. From \eqref{eq:Tk}, $T_k=\left\lfloor\frac{z_0-k\theta}{\theta}\right\rfloor=A_m-k$ for all $0\leq k\leq A_m$. From \eqref{eq:BS},
\[
\dd(m,T)=1+A_m+\sum_{k=0}^{A_m}\left(\left\lfloor\frac{y_0+kh}{y}\right\rfloor+A_m-2k\right).
\]
Setting $y_0=\overline{y_0}y+\widehat{y_0}$ with $0\leq\widehat{y_0}<y$ and $h=\overline{h}y+\widehat{h}$ with $0\leq\widehat{h}<y$, the above expression of $\dd(m,T)$ turns to be
\begin{align*}
\dd(m,T)&=1+A_m+\sum_{k=0}^{A_m}\left(\left\lfloor\frac{\widehat{y_0}+k\widehat{h}}{y}\right\rfloor+\overline{y_0}+k(\overline{h}-2)+A_m\right)\\
&=(1+A_m)(1+A_m+\overline{y_0})+(\overline{h}-2)\frac{A_m(1+A_m)}2+\sum_{k=0}^{A_m}\left\lfloor\frac{\widehat{y_0}+k\widehat{h}}{y}\right\rfloor.
\end{align*}
\item[(i.2)] Assume now $w>0$. Then,
\[
T_k=\min\left\{\left\lfloor\frac{x_0+k(l-w)}{w}\right\rfloor,\left\lfloor\frac{z_0-k\theta}{\theta}\right\rfloor\right\}=\min\left\{\left\lfloor\frac{x_0+kl}{w}\right\rfloor,\left\lfloor\frac{z_0}{\theta}\right\rfloor\right\}-k.
\]
The inequality $\left\lfloor\frac{x_0+kl}{w}\right\rfloor\geq\left\lfloor\frac{z_0}{\theta}\right\rfloor$ holds when either $\frac{x_0+kl}{w}\geq\frac{z_0}{\theta}$ or $n\leq\frac{x_0+kl}{w}<\frac{z_0}{\theta}<n+1$ for some $n\in\Nat$. The former holds when $k\geq k_0=\left\lceil\frac{z_0w-x_0\theta}{b}\right\rceil$, the latter holds whenever $0<\frac{z_0}{\theta}-\frac{x_0+kl}{w}<1\Leftrightarrow\frac{z_0w-x_0\theta}{b}-\frac wl<k<\frac{z_0w-x_0\theta}{b}$ (and, in this case, $\left\lfloor\frac{x_0+kl}{w}\right\rfloor=\left\lfloor\frac{z_0}{\theta}\right\rfloor$ holds). Notice that, from $0<\frac wl<1$, if there exists some $k_1\in\Nat$ such that $\frac{z_0w-x_0\theta}{b}-\frac wl<k_1<\frac{z_0w-x_0\theta}{b}$, this value $k_1$ must be unique and equality $\left\lfloor\frac{x_0+k_1l}{w}\right\rfloor=\left\lfloor\frac{z_0}{\theta}\right\rfloor$ holds. Thus, it follows that
\begin{equation}
T_k=\begin{cases}
\left\lfloor\frac{x_0+kl}{w}\right\rfloor-k&\text{if }k<k_0,\\[2mm]
A_m-k&\text{if }k\geq k_0.\label{eq:Tk2}
\end{cases}
\end{equation}
By Remark~\ref{rem:k0} we have $k_0\geq0$ and we consider three possible options.
\begin{itemize}
\item[(i.2.1)] Assume $k_0=0$. Then, $k\geq k_0$ for all $k$ and $T_k=A_m-k$ for all $k$. Therefore, the expression of $T_k$ is the same as in the previous case for all $k$. Thus the denumerant has the same expression as the previous case.

\item[(i.2.2)] Assume $1\leq k_0\leq A_m$. Now, the expression of $T_k$ changes upon the value of $k<k_0$ and $k\geq k_0$ according to \eqref{eq:Tk2}. Thus,
\begin{align*}
\dd(m,T)&=1+A_m+\sum_{k=0}^{A_m}\left(\left\lfloor\frac{y_0+kh}{y}\right\rfloor-k\right)+\sum_{k=0}^{k_0-1}\left(\left\lfloor\frac{x_0+kl}{w}\right\rfloor-k\right)+\sum_{k=k_0}^{A_m}\left(A_m-k\right)\\
&=1+A_m+\sum_{k=0}^{A_m}\left\lfloor\frac{\widehat{y_0}+k\widehat{h}}{y}\right\rfloor+\sum_{k=0}^{A_m}(\overline{y_0}+k(\overline{h}-2))+\sum_{k=0}^{k_0-1}\left\lfloor\frac{\widehat{x_0}+k\widehat{l}}{w}\right\rfloor\\
&\quad+\sum_{k=0}^{k_0-1}(\overline{x_0}+k\overline{l})+\sum_{k=k_0}^{A_m}A_m\\
&=(1+A_m)(1+A_m+\overline{y_0})+k_0(\overline{x_0}-A_m)+(\overline{h}-2)\frac{A_m(1+A_m)}2\nonumber\\
&\quad+\overline{l}\frac{(k_0-1)k_0}2+\sum_{k=0}^{A_m}\left\lfloor\frac{\widehat{y_0}+k\widehat{h}}y\right\rfloor+\sum_{k=0}^{k_0-1}\left\lfloor\frac{\widehat{x_0}+k\widehat{l}}w\right\rfloor,
\end{align*}
where $\overline{y_0},\widehat{y_0},\overline{h},\widehat{h},\overline{x_0},\widehat{x_0},\overline{l}$ and $\widehat{l}$ are those parameters defined in the statement (i.2.2) of the theorem.

\item[(i.2.3)] Assume $k_0>A_m$. Now, following \eqref{eq:Tk2}, we have $T_k=\left\lfloor\frac{x_0+kl}{w}\right\rfloor-k$ for all $k$. Then,
\begin{align*}
\dd(m,T)&=1+A_m+\sum_{k=0}^{A_m}\left(\left\lfloor\frac{y_0+kh}{y}\right\rfloor-k+\left\lfloor\frac{x_0+kl}{w}\right\rfloor-k\right)\\
&=1+A_m+\sum_{k=0}^{A_m}(\overline{y_0}+k\overline{h}+\overline{x_0}+k\overline{l}-2k)+\sum_{k=0}^{A_m}\left\lfloor\frac{\widehat{y_0}+k\widehat{h}}y\right\rfloor+\sum_{k=0}^{A_m}\left\lfloor\frac{\widehat{x_0}+k\widehat{l}}w\right\rfloor\\
&=(1+A_m)(1+\overline{x_0}+\overline{y_0})+(\overline{l}+\overline{h}-2)\frac{A_m(1+A_m)}2\\
&\quad+\sum_{k=0}^{A_m}\left\lfloor\frac{\widehat{x_0}+k\widehat{l}}w\right\rfloor+\sum_{k=0}^{A_m}\left\lfloor\frac{\widehat{y_0}+k\widehat{h}}y\right\rfloor,
\end{align*}
where $\overline{y_0},\widehat{y_0},\overline{h},\widehat{h},\overline{x_0},\widehat{x_0},\overline{l}$ and $\widehat{l}$ are the same as those defined in (i.2.2).\hfill$\square$
\end{itemize}
\end{itemize}

Notice how all the expressions of denumerant given by Theorem~\ref{teo:casi} contain sums of type $S^\pm$.

\subsection{Case (ii) $\theta=0$}\label{sec:casii}

This case is similar to the case (i). Now we have $w\neq0$, $\delta\neq0$, $A_m=\left\lfloor\frac{z_0}{\delta}\right\rfloor$ and $T_k=\left\lfloor\frac{x_0+kl}{w}\right\rfloor-k$ for all $k$. As in the previous case, we consider the following parameters
\begin{align}
x_0&=\overline{x_0}w+\widehat{x_0},\quad0\leq\widehat{x_0}<w,\nonumber\\
l&=\overline{l}w+\widehat{l},\qquad\,\,0\leq\widehat{l}<w,\nonumber\\
y_0&=\overline{y_0}y+\widehat{y_0},\quad\,\,0\leq\widehat{y_0}<y,\label{eq:paramii}\\
h&=\overline{h}y+\widehat{h},\qquad0\leq\widehat{h}<y.\nonumber
\end{align}
Defining $k_1=\left\lceil\frac{z_0y-y_0\delta}{a}\right\rceil$ (recall that $a=h\delta$) and using similar arguments of (i.2) in the proof of Theorem~\ref{teo:casi}, we have $S_k$ in \eqref{eq:Sk} turns to be
\begin{equation}
S_k=\begin{cases}
\left\lfloor\frac{y_0+kh}{y}\right\rfloor-k&\text{if }k<k_1,\\[2mm]
A_m-k&\text{if }k\geq k_1.\label{eq:Sk2}
\end{cases}
\end{equation}

The following result can be obtained using similar arguments like in the proof of Theorem~\ref{teo:casi}.

\begin{teo}\label{teo:casii}
Let us assume $\theta=0$. Set $A_m=\left\lfloor\frac{z_0}{\delta}\right\rfloor$. Then,
\begin{itemize}
\item[(ii.1)] if $y=0$, then
\begin{equation}
\dd(m,T)=(1+A_m)(1+A_m+\overline{x_0})+(\overline{l}-2)\frac{A_m(1+A_m)}2+\sum_{k=0}^{A_m}\left\lfloor\frac{\widehat{x_0}+k\widehat{l}}{w}\right\rfloor,\label{eq:casii1}
\end{equation}

\item[(ii.2)] if $y>0$, set $k_1=\left\lceil\frac{z_0y-y_0\delta}{a}\right\rceil$. Then,
\begin{itemize}
\item[(ii.2.1)] if $k_1=0$, then $\dd(m,T)$ has the same expression as in \eqref{eq:casii1}.
\item[(ii.2.2)] if $1\leq k_1\leq A_m$,
\begin{align}
\dd(m,T)&=(1+A_m)(1+A_m+\overline{x_0})+k_1(\overline{y_0}-A_m)+(\overline{l}-2)\frac{A_m(1+A_m)}2\nonumber\\
&\quad+\overline{h}\frac{(k_1-1)k_1}2+\sum_{k=0}^{A_m}\left\lfloor\frac{\widehat{x_0}+k\widehat{l}}w\right\rfloor+\sum_{k=0}^{k_1-1}\left\lfloor\frac{\widehat{y_0}+k\widehat{h}}y\right\rfloor,\label{eq:casii22}
\end{align}

\item[(ii.2.3)] if $k_1>A_m$, then the denumerant has the same expression as \eqref{eq:casi23}.
\end{itemize}
\end{itemize}
\end{teo}

\begin{rem}\sl
Although some expressions of Theorem~\ref{teo:casii} seem to be the same as some expressions of Theorem~\ref{teo:casi}, the value $A_m$ is not the same. In the former case we have $A_m=\left\lfloor\frac{z_0}{\theta}\right\rfloor$ and $A_m=\left\lfloor\frac{z_0}{\delta}\right\rfloor$ in the latter.
\end{rem}

\begin{rem}\sl\label{rem:k1}
In Theorem~\ref{teo:casii} we have $k_1\geq0$. Indeed, $k_1\geq0\Leftrightarrow-1<\frac{z_0y-y_0\delta}{h\delta}$ (recall that $h\delta=a$) and $-1<\frac{z_0y-y_0\delta}{h\delta}\Leftrightarrow0<z_0y+\delta(h-y_0)$. The factorization $(x_0,y_0,z_0)$ is the basic one of $m$ with respect to the L--shape $\HH$. Thus, $\sq{x_0,y_0}\in\HH\Rightarrow h>y_0$.
\end{rem}

\subsection{Case (iii) $\delta\theta>0$}
Now we have $A_m=\left\lfloor\frac{z_0}{\delta+\theta}\right\rfloor$, $a=h\delta+y\theta$ and $b=w\delta+l\theta$. There are four different options that give different expressions of $A=\sum_{k=0}^{A_m}S_k$ and $B=\sum_{k=0}^{A_m}T_k$ in \eqref{eq:BS},
\begin{itemize}
\item[](iii.1) $w=y=0$,
\item[](iii.2) $w\neq0$ and $y=0$,
\item[](iii.3) $w=0$ and $y\neq0$,
\item[](iii.4) $wy\neq0$.
\end{itemize}
Here we use the same notation as in \eqref{eq:paramii} plus the following one
\begin{align}
z_0&=\overline{z_{0,1}}\delta+\widehat{z_{0,1}},\quad 0\leq\widehat{z_{0,1}}<\delta,\nonumber\\
z_0&=\overline{z_{0,2}}\theta+\widehat{z_{0,2}},\quad 0\leq\widehat{z_{0,2}}<\theta,\nonumber\\
\delta&=\overline{\delta}\theta+\widehat{\delta},\quad\qquad 0\leq\widehat{\delta}<\theta,\label{eq:paramiii}\\
\theta&=\overline{\theta}\delta+\widehat{\theta},\quad\qquad 0\leq\widehat{\theta}<\delta,\nonumber
\end{align}

\begin{teo}\label{teo:casiii}
Let us assume the numerical $3$--semigroup $T=\sg{a,b,c}$ has related the L--shape $\HH=\ele(l,h,w,y)$. Consider $m=x_0a+y_0b+z_0c$, where $(x_0,y_0,z_0)$ is the basic factorization of $m$ with respect to $\HH$ in $T$. Then,
\begin{itemize}
\item[] (iii.1) if $w=y=0$,
\begin{align}
\dd(m,T)&=(1+A_m)(1+\overline{z_{0,1}}+\overline{z_{0,2}})-(\overline{\delta}+\overline{\theta}+2)\frac{A_m(1+A_m)}2\nonumber\\
&\quad+\sum_{k=0}^{A_m}\left\lfloor\frac{\widehat{z_{0,1}}-k\widehat{\theta}}{\delta}\right\rfloor+\sum_{k=0}^{A_m}\left\lfloor\frac{\widehat{z_{0,2}}-k\widehat{\delta}}{\theta}\right\rfloor,\label{eq:casiii1}
\end{align}

\item[] (iii.2) if $w\neq0$ and $y=0$, set $k_0=\left\lceil\frac{z_0w-x_0\theta}{b}\right\rceil$; then,
\begin{itemize}
\item[] (iii.2.1) if $k_0=0$, the denumerant has the same expression as in \eqref{eq:casiii1}. Otherwise, when $k_0>0$, we have
\begin{equation}
\dd(m,T)=(1+A_m)(1+\overline{z_{0,1}})-(1+\overline{\theta})\frac{A_m(1+A_m)}2+\sum_{k=0}^{A_m}\left\lfloor\frac{\widehat{z_{0,1}}-k\widehat{\theta}}{\delta}\right\rfloor+B\label{eq:casiii2B}
\end{equation}
where $B$ is defined through the rules of (iii.2.2) or (iii.2.3).
\item[] (iii.2.2) if $1\leq k_0\leq A_m$, we have
\begin{align}
B&=(1+A_m)\overline{z_{0,2}}+k_0(\overline{x_0}-\overline{z_{0,2}})+(\overline{l}+\overline{\delta})\frac{(k_0-1)k_0}2-(\overline{\delta}+1)\frac{A_m(1+A_m)}2\nonumber\\
&\quad+\sum_{k=0}^{k_0-1}\left\lfloor\frac{\widehat{x_0}+k\widehat{l}}{w}\right\rfloor+\sum_{k=k_0}^{A_m}\left\lfloor\frac{\widehat{z_{0,2}}-k\widehat{\delta}}{\theta}\right\rfloor\label{eq:casiii22}
\end{align}
\item[] (iii.2.3) if $k_0>A_m$, then
\begin{equation}
B=(1+A_m)\overline{x_0}+(\overline{l}-1)\frac{A_m(1+A_m)}2+\sum_{k=0}^{A_m}\left\lfloor\frac{\widehat{x_0}+k\widehat{l}}{w}\right\rfloor,\label{eq:casiii23}
\end{equation}
\end{itemize}

\item[] (iii.3) if $w=0$ and $y\neq0$, define $k_1=\left\lceil\frac{z_0y-y_0\delta}{a}\right\rceil$; then,
\begin{itemize}
\item[] (iii.3.1) if $k_1=0$, the denumerant has the same expression as in \eqref{eq:casiii1}.

Otherwise, when $k_1>0$, we have
\begin{equation}
\dd(m,T)=(1+A_m)(1+\overline{z_{0,2}})-(\overline{\delta}+1)\frac{A_m(1+A_m)}2+\sum_{k=0}^{A_m}\left\lfloor\frac{\widehat{z_{0,2}}-k\widehat{\delta}}{\theta}\right\rfloor+A\label{eq:casiii3A}
\end{equation}
where $A$ is defined by (iii.3.2) or (iii.3.3).
\item[] (iii.3.2) if $1\leq k_1\leq A_m$, then
\begin{align}
A&=(1+A_m)\overline{z_{0,1}}+k_1(\overline{y_0}-\overline{z_{0,1}})+(\overline{\theta}+\overline{h})\frac{(k_1-1)k_1}2-(\overline{\theta}+1)\frac{A_m(1+A_m)}2\nonumber\\
&\quad+\sum_{k=0}^{k_1-1}\left\lfloor\frac{\widehat{y_0}+k\widehat{h}}{y}\right\rfloor+\sum_{k=k_1}^{A_m}\left\lfloor\frac{\widehat{z_{0,1}}-k\widehat{\theta}}{\delta}\right\rfloor
\end{align}

\item[] (iii.3.3) if $k_1>A_m$, then
\begin{equation}
A=(1+A_m)\overline{y_0}+(\overline{h}-1)\frac{A_m(1+A_m)}2+\sum_{k=0}^{A_m}\left\lfloor\frac{\widehat{y_0}+k\widehat{h}}{y}\right\rfloor,\label{eq:casiii33}
\end{equation}
\end{itemize}

\item[] (iii.4) if $wy\neq0$, define $k_0$ and $k_1$ as in (iii.2) and (iii.3), respectively; then
\begin{equation}
\dd(m,T)=1+A_m+A+B,
\end{equation}
where $A$ and $B$ are ruled by the following expressions, depending on $k_0$ and $k_1$.
\begin{itemize}
\item[-] If $k_1=0$, then
\begin{equation}
A=(1+A_m)\overline{z_{0,1}}-(\overline{\theta}+1)\frac{A_m(1+A_m)}2+\sum_{k=0}^{A_m}\left\lfloor\frac{\widehat{z_{0,1}}-k\widehat{\theta}}{\delta}\right\rfloor.\label{eq:casiii4A1}
\end{equation}
\item[-] If $1\leq k_1\leq A_m$, then
\begin{align}
A&=(1+A_m)\overline{z_{0,1}}-(\overline{\theta}+1)\frac{A_m(1+A_m)}2+k_1(\overline{y_0}-\overline{z_{0,1}})+(\overline{h}+\overline{\theta})\frac{(k_1-1)k_1}2\nonumber\\
&\quad+\sum_{k=0}^{k_1-1}\left\lfloor\frac{\widehat{y_0}+k\widehat{h}}{y}\right\rfloor+\sum_{k=k_1}^{A_m}\left\lfloor\frac{\widehat{z_{0,1}}-k\widehat{\theta}}{\delta}\right\rfloor.\label{eq:casiii4A2}
\end{align}
\item[-] If $k_1>A_m$, then $A$ has the same expression as in \eqref{eq:casiii33}.
\item[-] If $k_0=0$, then
\begin{equation}
B=(1+A_m)\overline{z_{0,2}}-(\overline{\delta}+1)\frac{A_m(1+A_m)}2+\sum_{k=0}^{A_m}\left\lfloor\frac{\widehat{z_{0,2}}-k\widehat{\delta}}{\theta}\right\rfloor.\label{eq:casiii4B1}
\end{equation}
\item[-] If $1\leq k_0\leq A_m$, then $B$ has the same expression as that in \eqref{eq:casiii22}.
%\begin{align}
%B&=(1+A_m)\overline{z_{0,2}}-(\overline{\delta}+1)\frac{A_m(1+A_m)}2+k_0(\overline{x_0}-\overline{z_{0,2}})+(\overline{\delta}+\overline{l})\frac{(k_0-1)k_0}2\nonumber\\
%&\quad+\sum_{k=0}^{k_0-1}\left\lfloor\frac{\widehat{x_0}+k\widehat{l}}{w}\right\rfloor+\sum_{k=k_0}^{A_m}\left\lfloor\frac{\widehat{z_{0,2}}-k\widehat{\delta}}{\theta}\right\rfloor.\label{eq:casiii4B2}
%\end{align}
\item[-] If $k_0>A_m$, then $B$ has the same expression as in \eqref{eq:casiii23}.
\end{itemize}
\end{itemize}
\end{teo}
\noindent{\bf Proof}: For the stated values of $k_0$ and $k_1$ (recall that now we have $a=h\delta+y\theta$ and $b=w\delta+l\theta$), from \eqref{eq:Sk} and \eqref{eq:Tk}, we still have
\[
S_k=\begin{cases}
\left\lfloor\frac{y_0+kh}{y}\right\rfloor-k&\text{if }k<k_1,\\[2mm]
A_m-k&\text{if }k\geq k_1,
\end{cases}
\quad\text{ and }\quad
T_k=\begin{cases}
\left\lfloor\frac{x_0+kl}{w}\right\rfloor-k&\text{if }k<k_0,\\[2mm]
A_m-k&\text{if }k\geq k_0.
\end{cases}
\]
Then, all the expressions of the statement are obtained using the same arguments of the proof of Theorem~\ref{teo:casi}.\hfill$\square$

\begin{rem}\sl
Similar arguments of remarks \ref{rem:k0} and \ref{rem:k1} leave to $k_0\geq0$ and $k_1\geq0$.
\end{rem}

\begin{rem}\sl\label{rem:suma}
Although in the statement of Theorem~\ref{teo:casiii} appear sums like $s=\sum_{k=k_0}^{A_m}\left\lfloor\frac{\widehat{z_{0,2}}-k\widehat{\delta}}{\theta}\right\rfloor$ that is not of type $S^\pm$ (the sum do not begin at $k=0$), we can reduce it to one sum of type $S^\pm$. Indeed, taking a generic sum $\sum_{k=n_1}^{n_2}\left\lfloor\frac{s\pm kt}{q}\right\rfloor$ with $0\leq s,t<q$, and changing the summation index, $u=k-n_1$, we obtain an $S^\pm$ sum
\[
\sum_{u=0}^{n_2-n_1}\left\lfloor\frac{s\pm n_1t\pm ut}{q}\right\rfloor=\overline{\alpha}(1+n_2-n_1)+\sum_{u=0}^{n_2-n_1}\left\lfloor\frac{\widehat{\alpha}\pm ut}{q}\right\rfloor
\]
with $\alpha=s\pm n_1t$, $\alpha=\overline{\alpha}q+\widehat{\alpha}$ and $0\leq\widehat{\alpha}<q$.
\end{rem}

\section{Discrete sums $S^\pm$}\label{sec:Spm}

Let us denote the {\em discrete sum $S^\pm$} by
\begin{equation}
S^{\pm}(s,t,q,N)=\sum_{k=0}^N\left\lfloor\frac{s\pm kt}q\right\rfloor,\qquad0\leq s,t<q.\label{eq:DS}
\end{equation}
These type of sums appear to be a main tool for computing denumerants, as it has been seen in the previous section. In this section we study some properties of $S^\pm$ in order to obtain an efficient numerical calculation of it. This calculation will be done in a {\em discrete Lebesgue--like} sense.

\subsection{$S^+$ sums}
Consider the function $f(x)=\left\lfloor\frac{s+xt}{q}\right\rfloor$ that defines the general term of an $S^+(s,t,q,N)$ sum.

\begin{defi}
Let us define the $k$--interval $I_k\subset[0,N]$ by $I_k=\{x\in[0,N]:~ f(x)=k\}$. A $k$--interval $I_k$ is called an hS-type interval if $|I_k\cap\Nat|=\lceil\frac qt\rceil$.
\end{defi}

\begin{lem}\label{lem:intk-basic}
Given a $k$--interval $I_k$, we have
\begin{itemize}
\item[(i)] $I_k=[x_k,x_{k+1})$ with $x_k=\frac{kq-s}t$.
\item[(ii)] $\lfloor\frac qt\rfloor\leq|I_k\cap\Nat|\leq\lceil\frac qt\rceil$ holds except, eventually, the first and/or last intervals.
\end{itemize}
\end{lem}
\noindent\textbf{Proof}: Item (i) comes directly from the expression of $f$. A real interval, $I=[\alpha,\beta)$ of length $\ell=\beta-\alpha$, contains at least $\lfloor\ell\rfloor$ integers and it contains at most $\lceil\ell\rceil$ integers. Item (ii) comes from the length of $|I_k|=x_{k+1}-x_k=\frac qt$. $\square$

We discuss the value of $S^+(s,t,q,N)$ depending on the following three subcases
\begin{itemize}
\item[(a)] $t\mid q$,
\item[(b)] $t\nmid q$ and $\gcd(t,q)=1$,
\item[(c)] $t\nmid q$ and $\gcd(t,q)=g>1$.
\end{itemize}

\subsubsection{Assume $t\mid q$}

The maximum value attained by $f$ in $[0,N]$ is
\begin{equation}
M=\left\lfloor\frac{s+Nt}{q}\right\rfloor\quad\text{at}\quad x_M=\frac{Mq-s}t.\label{eq:MxM-1}
\end{equation}

\begin{figure}[h]
\centering
\includegraphics[width=0.9\linewidth]{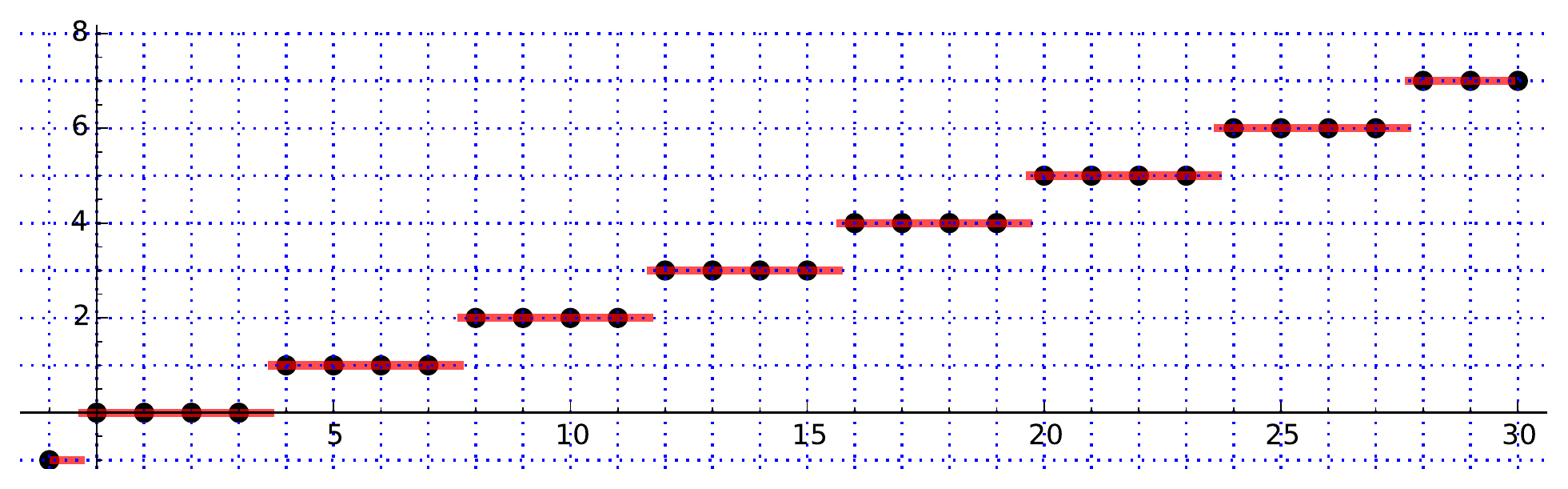}
\caption{$k$-intervals for $s=1$, $t=3$, $q=12$ and $N=30$}
\label{fig:s1t3q12N30}
\end{figure}

\begin{exa}\small\sl
Figure~\ref{fig:s1t3q12N30} shows the case $s=1$, $t=3$, $q=12$ and $N=30$. All $k$-intervals are hS-type ones. The distribution of integrals values is the same in each interval.
\end{exa}

\begin{teo}\label{teo:sm-1}
Let assume $t\mid q$. Then,
\[
S^+(s,t,q,N)=\frac qt\;\frac{M(M-1)}2+M(N-\lceil x_M\rceil+1).
\]
\end{teo}
\noindent\textbf{Proof}: Let us denote $q=t\overline{q}$. By Lemma~\ref{lem:intk-basic}, each $I_k$ interval is a hS--type interval, that is $|I_k|=\overline{q}$ (except, perhaps, the first $I_0$ and the last one $I_M$). We divide the interval $I=[0,N]$ in three regions $I=I_0\cup J\cup I_M$, where $M$ is the maximum value attained by the function $f$ in $I$.

Let us denote $n=\frac{x_M-x_1}{\overline{q}}$. Then, there are $n$ $k$-intervals different from $I_0$ and $I_M$ in $I$, i.e. $I_1=[x_1,x_2),\ldots,I_{M-1}=[x_{M-1},x_M)$ (here $n=M-1$ holds). The last interval $I_M$ can be eventually one point (that is $I_M=\{N\}$). The sum is
\[
S^+(s,t,q,N)=\sum_{k=0}^N f(k)=0+\sum_{k=x_1}^{x_{M}-1}f(k)+\sum_{k=x_M}^N f(k).
\]
Now we add these values like a discrete Lebesgue--like sum
\[
\sum_{k=x_1}^{x_{M}-1}f(k)=\sum_{j=1}^{M-1}|I_j\cap\Nat|j=\sum_{j=1}^{M-1}\overline{q}j=\overline{q}\;\frac{(M-1)M}2.
\]
Finally, we have to add $\sum_{k=x_M}^N f(k)=|I_M\cap\Nat|M$. The number of integral points in $I_M$ is $|I_M\cap\Nat|=N-\lceil x_M\rceil+1$. Thus, the value of $S^+(s,t,q,N)$ is the stated one. $\square$

\subsubsection{Assume $t\nmid q$ and $\gcd(t,q)=1$}

Assume $t\nmid q$. We use the notation
\begin{align}
q&=\overline{q}t+\hat{q},\quad 1\leq\hat{q}<t,\label{eq:pars1}\\
s&=\overline{s}t+\hat{s},\quad 0\leq\hat{s}<t,\label{eq:pars2}\\
S(s,t,q,N)&=\overline{q}\frac{(M-1)M}2+M(N-\lceil x_M\rceil+1).\label{eq:Smes}
\end{align}

\begin{defi}\label{defi:subhSind}
Given a set $A\subset\Nat$, a subset $J\subset A$ of hS indices is a set of ordered indices in $A$ of hS--type intervals. We define $S_J=\sum_{k\in J}k$.
\end{defi}

\begin{lem}\label{lem:hS-car}
Assume $t\nmid q$. Then, $I_k\subset[0,N]$ is an hS--type interval if and only if
\[
(\hat{s}-k\hat{q})\mmod{t}<\hat{q}.
\]
\end{lem}
\noindent\textbf{Proof}: The modulo in the statement is taken from the set of residues $\{0,1,\ldots,t-1\}$. Notice that $|I_k|=x_{k+1}-x_k=\overline{q}+\frac{\hat{q}}t$. The interval $I_k$ is hS--type if and only if $\lceil x_k\rceil-x_k<\frac{\hat{q}}t$ (so, the maximum number of integral values are located in $I_k$). This condition can be restated in a more numerically stable relation. From
\[
x_k=\frac{kq-s}t=k\overline{q}-\overline{s}-\frac{\hat{s}-k\hat{q}}t,
\]
putting $\hat{s}-k\hat{q}=\alpha t+\beta$ with $0\leq\beta<t$, we have $x_k=n-\frac{\beta}t$. Thus, inequality $\lceil x_k\rceil-x_k<\frac{\hat{q}}t$ holds if and only if $\beta<\hat{q}$. Equivalently $(\hat{s}-k\hat{q})\mmod{t}<\hat{q}$. $\square$

\begin{figure}[h]
\centering
\includegraphics[width=0.9\linewidth]{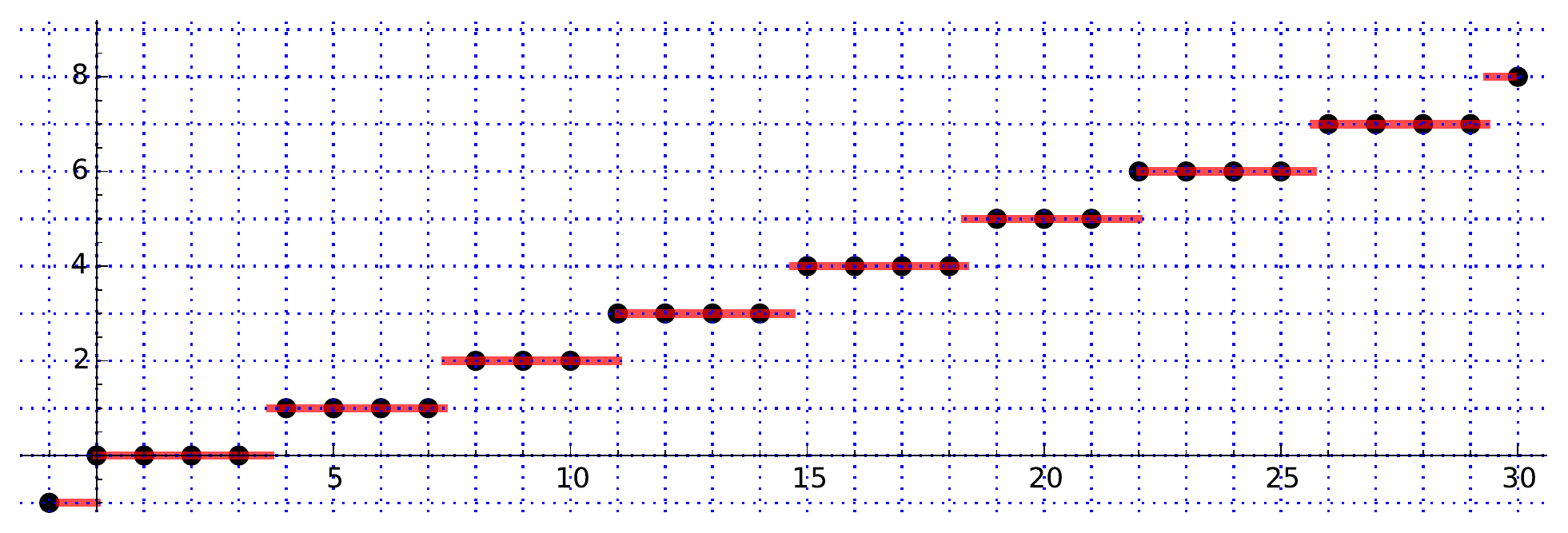}
\caption{$k$-intervals for $s=0$, $t=3$, $q=11$ and $N=30$}
\label{fig:kintervals1}
\end{figure}

\begin{exa}\small\sl
Let us consider $s=0$, $t=3$, $q=11$ and $N=30$. Figure~\ref{fig:kintervals1} shows all $k$-intervals in this case. Notice the hS-type intervals for $k\in\{0,1,3,4,6,7\}$ (when $0\leq k\leq M=8$ is a solution of $-2k\mmod{3}<2$). Thus, only the intervals $I_2$ and $I_5$ are not of type hS in Figure~\ref{fig:kintervals1}.
\end{exa}

The distribution pattern of integral values inside the $k$-intervals is also ruled modulo $t$. This fact is detailed in the following result.

\begin{lem}\label{lem:Ikper}
Let $I_k$ and $I_{k+T}$, $T>0$, be two intervals with the same distribution of integral values. Then,
\begin{itemize}
\item[(i)] $T\equiv0\mmod{t}$.
\item[(ii)] The minimum value of $T$ is $t$.
\end{itemize}
\end{lem}
\noindent\textbf{Proof}: In particular, $\lceil x_{k+T}\rceil-x_{k+T}=\lceil x_{k}\rceil-x_{k}$ holds. Then, $x_{k+T}-x_{k}\in\Z$, that is (recalling that $\gcd(q,t)=1$)
\[
x_{k+T}-x_{k}=T\;\frac qt\in\Z\Leftrightarrow T\equiv0\mmod{t}.
\]
The minimum $T>0$ for the value $T\frac qt$ to be an integer is $T=t$. $\square$

\begin{cor}\label{cor:tperiod}
The distribution of integral values in the $k$-intervals has period $t$.
\end{cor}

Although Lemma~\ref{lem:hS-car} and Lemma~\ref{lem:Ikper} give a characterization of hS--intervals, we need a more accurate description of these intervals. This description will be used to efficiently obtain a subset of hS--indices $J$ of Definition~\ref{defi:subhSind}. Indeed, from Lemma~\ref{lem:hS-car}, the set $J$ can be parameterized by
\begin{equation}
J=\{\widehat{q}^{-1}(\widehat{s}-i)\mmod{t}|~0\leq i\leq\widehat{q}-1\}.\label{eq:Jpnosorted}
\end{equation}
Notice that $J\neq\emptyset$ because $\widehat{q}\geq1$ ($t\nmid q$ and $\gcd(t,q)=1$) and $\widehat{q}^{-1}\widehat{s}\in J$ always. This parameterization is useless, from the point of view of numerical efficiency, whenever we need the elements of $J$ to be sorted. Noting that elements of $J$ are sorted by a rule defined by two moduli, $\widehat{q}$ and $t$, we can obtain a sorted parameterization of $J$. For instance (notice that $\gcd(t,\widehat{q})=1$)
\begin{equation}
J=\{[\widehat{q}^{-1}(\widehat{s}-(\widehat{s}+iu)\mmod{\widehat{q}})]\mmod{t}|~0\leq i\leq\widehat{q}-1\},\qquad u\equiv t\mmod{\widehat{q}}\label{eq:Jpsorted}
\end{equation}
is an example of such parameterization.

\begin{teo}\label{teo:sm-2}
Assume $t\nmid q$ and $\gcd(t,q)=1$. Consider the set of hS-type indices $J=\{j_0,\ldots,j_{\widehat{q}-1}\}\subset\{0,\ldots,t-1\}$ and $S(s,t,q,N)$ given by expression \eqref{eq:Smes}. Then,
\begin{itemize}
\item[(a)] If $j_0\geq M$, then $S^+(s,t,q,N)=S(s,t,q,N)$ holds.
\item[(b)] If $j_0<M$, there are three different cases.
\begin{itemize}
\item[(b.1)] If $j_{\widehat{q}-1}\geq M$, consider the subset of hS indices $K\subset\{0,\ldots,M-1\}$. Then,  $S^+(s,t,q,N)=S(s,t,q,N)+S_K$ holds.
\item[(b.2)] If $j_{\widehat{q}-1}<M$ and $j_0+t\geq M$, then $S^+(s,t,q,N)=S(s,t,q,N)+S_J$ holds.
\item[(b.3)] If $j_{\widehat{q}-1}<M$ and $j_0+t<M$, set $u=\left\lfloor\frac{M-1}t\right\rfloor$ and consider the set of hS-type indices $K\subset\{j_0+ut,\ldots,M-1\}$. Then,
\[
S^+(s,t,q,N)=S(s,t,q,N)+uS_J+\widehat{q}t\frac{(u-1)u}2+S_K.
\]
\end{itemize}
\end{itemize}
\end{teo}
\noindent\textbf{Proof}: $S^+(s,t,q,N)$ can be calculated from $S$ plus all additional summands corresponding to hS-type intervals. That is, each hS-type interval $I_j$ has an additional value $j$ which must be added to $S$.

(a) When $j_0\geq M$, there is no hS-type interval in $[0,M)$. Thus, all $k$-intervals in this region has $\overline{q}$ integral values. Then, $S^+$ has the same expression as in Theorem~\ref{teo:sm-1} replacing $\frac qt$ by $\overline{q}$, that is $S^+(s,t,q,N)=S(s,t,q,N)$.

(b) When $j_0<M$ there are hS-type intervals in $[0,M)$. So, we also have to add all hS-type indices contained in $[0,M)$ for obtaining $S^+$.
\begin{itemize}
\item[(b.1)] Assume $j_{\widehat{q}-1}\geq M$. Consider the set of hS-type indices $K\subset\{0,\ldots,M-1\}$. There are no more hS-type indices to consider and $S^+(s,t,q,N)=S(s,t,q,N)+S_K$.

\item[(b.2)] Assume $j_{\widehat{q}-1}<M$ and $j_0+t\geq M$. Then, by Lemma~\ref{lem:hS-car}, all hS-type indices are $J$. Thus, $S^+(s,t,q,N)=S(s,t,q,N)+S_J$ holds.

\item[(b.3)] Assume $j_{\widehat{q}-1}<M$ and $j_0+t<M$. By Lemma~\ref{lem:hS-car}, the behaviour of the $k$-intervals is $t$-periodic. The maximum number of periods included in the set of indices $A=\{0,\ldots,M-1\}$ is $u=\left\lfloor\frac{M-1}t\right\rfloor$. That is, all the elements in $\{j_0,\ldots,j_{\widehat{q}-1},j_0+t,\ldots,j_{\widehat{q}-1}+t,\ldots,j_0+(u-1)t,\ldots,j_{\widehat{q}-1}+(u-1)t\}$ are hS-type indices. The remaining hS-type indices are located in the set of hS indices $K\subset\{ut,\ldots,M-1\}$. Therefore,
\begin{align*}
S^+(s,t,q,N)&=S+\sum_{l=0}^{u-1}\left(\sum_{j\in J}(lt+j)\right)+S_K=S+\sum_{l=0}^{u-1}\left(lt|J|+S_J\right)+S_K\\
&=S+|J|t\sum_{l=0}^{u-1}l+uS_J+S_K=S+|J|t\frac{(u-1)u}2+uS_J+S_K.
\end{align*}
The statement follows from the identity $|J|=\widehat{q}$. $\square$
\end{itemize}

\begin{rem}\label{rem:tperiod}\sl
The sets of hS-type indices $J$ and $K$ of Theorem~\ref{teo:sm-2} are obtained at time cost $O(\widehat{q})$, in the worst case. The first and last elements of $J$, $j_0$ and $j_{\widehat{q}-1}$, can be obtained at constant time cost from the (sorted) parameterization \eqref{eq:Jpsorted} of $J$.
\end{rem}

\subsubsection{Assume $t\nmid q$ and $\gcd(t,q)=g>1$}

When $t\nmid q$ and $\gcd(t,q)=g>1$, we have $x_{h+1}-x_k=\frac qt=\frac{\tilde{q}}{\tilde{t}}$, where $\tilde{t}=\frac tg$ and $\tilde{q}=\frac qg$.

\begin{lem}\label{lem:Ikper-2}
Assume $t\nmid q$ and $\gcd(t,q)=g>1$. Let's assume $I_k$ and $I_{k+T}$, $T>0$, are two intervals with the same distribution of integral values. Then,
\begin{itemize}
\item[(i)] $T\equiv0\mmod{\tilde{t}}$.
\item[(ii)] The minimum value of $T$ is $\tilde{t}$.
\end{itemize}
\end{lem}
\noindent\textbf{Proof}: This lemma follows from the proof of Lemma~\ref{lem:Ikper} with the additional identity $\frac qt=\frac{\tilde{q}}{\tilde{t}}$, $\gcd(\tilde{t},\tilde{q})=1$. $\square$ 

In particular, Lemma~\ref{lem:Ikper-2} ensures that the distribution of integrals values of $k$-intervals in $[0,N]$ has period $\tilde{t}$. Now, by Lemma~\ref{lem:Ikper-2}, detecting hS-type intervals is done as follows. Set
\begin{align}
s&=\overline{s}g+s_g,\quad\;\, 0\leq s_g<g,\nonumber\\
\overline{s}&=\overline{\overline{s}}\tilde{t}+\hat{\hat{s}},\qquad 0\leq\hat{\hat{s}}<\tilde{t},\label{eq:par2}\\
\tilde{q}&=\overline{\overline{q}}\tilde{t}+\hat{\hat{q}},\qquad 1\leq\hat{\hat{q}}<\tilde{t}.\nonumber
\end{align}
Then, $I_k$ is an hS interval if and only if
\begin{equation}
(\hat{\hat{s}}-k\hat{\hat{q}})\mmod{\tilde{t}}<\hat{\hat{q}},\label{eq:hS-car2}
\end{equation}
that is similar to the characterization given in Lemma~\ref{lem:hS-car}.

\begin{figure}[h]
\centering
\includegraphics[width=0.9\linewidth]{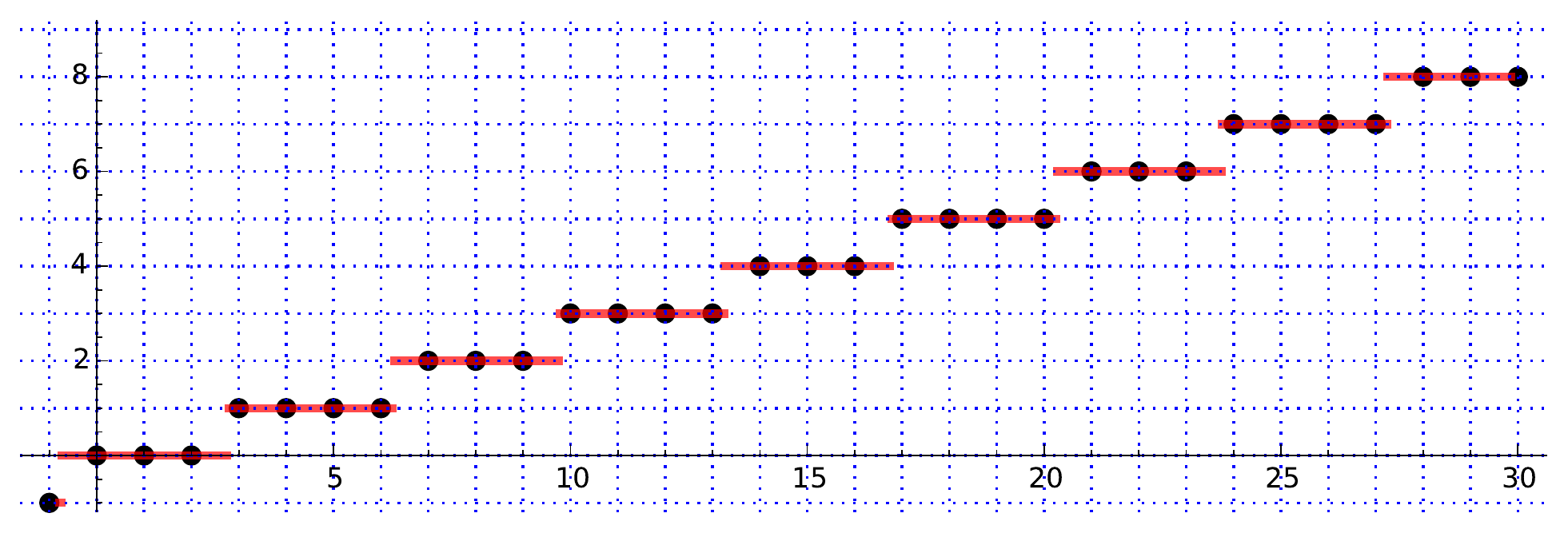}
\caption{$k$-intervals for $s=3$, $t=4$, $q=14$ and $N=30$}
\label{fig:s3t4q14N30}
\end{figure}

\begin{exa}\small\sl
Figure~\ref{fig:s3t4q14N30} shows $k$-intervals in the case $s=3$, $t=4$, $q=14$ and $N=30$. Here the distribution of integrals values in $k$-intervals has period $\tilde{t}=2$, that is hS-type intervals follow the rule \eqref{eq:hS-car2}, i.e. the set hS-type indices in $[0,8)$ is $\{1,3,5,7\}$.
\end{exa}

\begin{rem}\sl
Sorted and non sorted parameterization of $J$ can also be obtained as in \eqref{eq:Jpnosorted} and \eqref{eq:Jpsorted}. The same expressions hold  replacing $t$ by $\tilde{t}$, $\widehat{s}$ by $\hat{\hat{s}}$ and $\widehat{q}$ by $\hat{\hat{q}}$.
\end{rem}

Now, we denote
\begin{equation}\label{eq:s2}
S(s,t,q,N)=\left\lfloor\frac{\tilde{q}}{\tilde{t}}\right\rfloor\;\frac{M(M-1)}2+M(N-\lceil x_M\rceil+1)
\end{equation}
and similar results are obtained from the $\tilde{t}$--periodicity of the hS-type intervals ruled by \eqref{eq:hS-car2}.

\begin{teo}\label{teo:sm-3}
Assume $t\nmid q$ and $\gcd(t,q)=g>1$. Set $\tilde{t}=t/g$ and $\tilde{q}=q/g$. Then, interchanging $t$ by $\tilde{t}$ and $\widehat{q}$ by $\hat{\hat{q}}$, statements of Theorem~\ref{teo:sm-2} hold.
\end{teo}

\begin{rem}\label{rem:MxM}\sl
Notice that $M$ and $x_M$ are also calculated like in \eqref{eq:MxM-1}, i.e. using $t$ and $q$ (not $\tilde{t}$ and $\tilde{q}$). 
\end{rem}

\begin{rem}\label{rem:tcircperiod}\sl
Now, the sets of hS indices $J$ and $K$ are computed using \eqref{eq:hS-car2} at time cost $O(\tilde{t})$.
\end{rem}

\subsection{$S^-$ sums}

The {\em minus sums}
\begin{equation}
S^-(s,t,q,N)=\sum_{k=0}^N\left\lfloor\frac{s-kt}{q}\right\rfloor\label{eq:smenys}
\end{equation}
share some behaviour with plus sums $S^+(s,t,q,N)$. We can define by analogy $k$-intervals $I_k=(x_k,x_{k+1}]\subset\Real$ (those intervals such that $g(x)=\left\lfloor\frac{s-xt}{q}\right\rfloor=k$ for $x\in I_k$) with $x_k=\frac{s-(k+1)q}t$ and $x_{k+1}=\frac{s-kq}t$. hS-type intervals are also defined to be those $I_k$ with $|I_k\cap\Z|=\left\lceil\frac qt\right\rceil$. We denote now
\begin{equation}
M=-\left\lfloor\frac{s-Nt}q\right\rfloor\quad\text{and}\quad x_M=\frac{s+(M-1)q}t,\label{eq:MxM-2}
\end{equation}
that are the analog to \eqref{eq:MxM-1} for $S^+$. Also three cases are taken into account now, i.e. $t\mid q$, $t\nmid q$ with $\gcd(t,q)=1$ and $t\nmid q$ with $\gcd(t,q)=g>1$. We give here, without proof, the main results for computing $S^-$ sums.

%\begin{figure}[h]
%\centering
%\includegraphics[width=0.9\linewidth]{s3_t2_q8_N30.pdf}
%\caption{$k$-intervals for $s=3$, $t=2$, $q=8$ and $N=30$}
%\label{fig:s3t2q8N30}
%\end{figure}

When $t\mid q$, all intervals are hS-type ones and have the same distribution of integral values. The following result can be proved using similar arguments as in Theorem~\ref{teo:sm-1}.

\begin{teo}\label{teo:smenys-1}
Assume $t\mid q$. Then,
\[
S^-(s,t,q,N)=-\frac qt\;\frac{(M-1)M}2-M(N-\lfloor x_M\rfloor).
\]
\end{teo}

%\begin{figure}[h]
%\centering
%\includegraphics[width=0.9\linewidth]{s3_t3_q10_N30.pdf}
%\caption{$k$-intervals for $s=3$, $t=3$, $q=10$ and $N=30$}
%\label{fig:s3t2q8N30}
%\end{figure}

When $t\nmid q$ and $\gcd(t,q)=1$, we also use the notation $\widehat{q}$ and $\widehat{s}$ defined in \eqref{eq:pars1} and \eqref{eq:pars2}.

\begin{lem}\label{lem:hS-car-menys}
Assume $t\nmid q$. Then, $I_k\subset[0,N]$ is an hS--type interval if and only if
\[
(\hat{s}+k\hat{q})\mmod{t}<\hat{q}.
\]
\end{lem}

Lemma~\ref{lem:hS-car-menys} allows a non sorted parameterization of the set of hS-type indices $J\subset\{0,\ldots,t-1\}$, that is
\begin{equation}
J=\{\widehat{q}^{-1}(i-\widehat{s})\mmod{t}|~0\leq i\leq\widehat{q}-1\},\label{eq:Jmnosorted}
\end{equation}
which is an analogous expressions to \eqref{eq:Jpnosorted} for plus sums. A sorted parameterization of $J$ is given by (now $u\equiv-t\mmod{\widehat{q}}$)
\begin{equation}
J=\left\{\widehat{q}^{-1}[(\widehat{s}+iu)\mmod{\widehat{q}}-\widehat{s}]\mmod{t}|~0\leq i\leq\widehat{q}-1\right\}\quad\text{ if }\widehat{s}<\widehat{q}\label{eq:Jmsorted1}
\end{equation}
and
\begin{equation}
J=\{\widehat{q}^{-1}[(\widehat{s}+(i+1)u)\mmod{\widehat{q}}-\widehat{s}]\mmod{t}|~0\leq i\leq\widehat{q}-1\}\quad\text{ if }\widehat{s}\geq\widehat{q}.\label{eq:Jmsorted2}
\end{equation}
In any case, as it has been done before, the sorted elements of $J$ will be denoted by $J=\{j_0,\ldots,j_{\widehat{q}-1}\}$.

The distribution of integral values in $k$-intervals also has period $t$ on the indices $k$ like in the plus sums. Let us denote the sum
\begin{equation}
S(s,t,q,N)=-\overline{q}\;\frac{(M-1)M}2-M(N-\lfloor x_M\rfloor)\label{eq:Smenys}
\end{equation}
which corresponds to $S^-$ when there is no hS-type $k$--interval, similar to \eqref{eq:Smes} for $S^+$. The following result is the analog of Theorem~\ref{teo:sm-2} for $S^+$.

\begin{teo}\label{teo:smenys-2}
Assume $t\nmid q$ and $\gcd(t,q)=1$. Consider the set of hS-type indices $J=\{j_0,\ldots,j_{\widehat{q}-1}\}\subset\{0,\ldots,t-1\}$ and $S(s,t,q,N)$ given by expression \eqref{eq:Smenys}. Then,
\begin{itemize}
\item[(a)] If $j_0\geq M$, then $S^+(s,t,q,N)=S(s,t,q,N)$ holds.
\item[(b)] If $j_0<M$, there are three different cases:
\begin{itemize}
\item[(b.1)] If $j_{\widehat{q}-1}\geq M$, consider the subset of hS indices $K\subset\{0,\ldots,M-1\}$. Then,  $S^+(s,t,q,N)=S(s,t,q,N)-S_K$ holds.
\item[(b.2)] If $j_{\widehat{q}-1}<M$ and $j_0+t\geq M$, then $S^+(s,t,q,N)=S(s,t,q,N)-S_J$ holds.
\item[(b.3)] If $j_{\widehat{q}-1}<M$ and $j_0+t<M$, set $u=\left\lfloor\frac{M-1}t\right\rfloor$ and consider the set of hS-type indices $K\subset\{j_0+ut,\ldots,M-1\}$. Then,
\[
S^+(s,t,q,N)=S(s,t,q,N)-uS_J-\widehat{q}t\frac{(u-1)u}2-S_K.
\]
\end{itemize}
\end{itemize}
\end{teo}

When $t\nmid q$ and $\gcd(t,q)=g>1$, we denote $\tilde{t}=t/g$ and $\tilde{q}=q/g$. The value $\overline{q}$ in (\ref{eq:Smenys}) is the same. i.e. $\overline{q}=\lfloor\tilde{q}/\tilde{t}\rfloor=\lfloor q/t\rfloor$. Using the same notation as in \eqref{eq:par2}, the analog to Lemma~\ref{lem:hS-car-menys} is
\begin{equation}
(\hat{\hat{s}}+k\hat{\hat{q}})\mmod{\tilde{t}}<\hat{\hat{q}}\label{eq:hS-car-menys2}
\end{equation}
and non sorted and sorted characterizations of $J$, \eqref{eq:Jmnosorted},  \eqref{eq:Jmsorted1} and \eqref{eq:Jmsorted2}, have the same expressions by replacing $u$ by $u\equiv-\tilde{t}\mmod{\hat{\hat{q}}}$, $t$ by $\tilde{t}$,  $\widehat{s}$ by $\hat{\hat{s}}$ and $\widehat{q}$ by $\hat{\hat{q}}$.

\begin{teo}\label{teo:smenys-3}
Assume $t\nmid q$ and $\gcd(t,q)=g>1$. Set $\tilde{t}=t/g$ and $\tilde{q}=q/g$. Then, interchanging $t$ by $\tilde{t}$ and $\widehat{q}$ by $\hat{\hat{q}}$, statements of Theorem~\ref{teo:smenys-2} hold.
\end{teo}

%\begin{figure}[h]
%\centering
%\includegraphics[width=0.9\linewidth]{s1_t10_q14_N30.pdf}
%\caption{$k$-intervals for $s=1$, $t=10$, $q=14$ and $N=30$}
%\label{fig:s1t10q14N30}
%\end{figure}

Remarks \ref{rem:MxM} and \ref{rem:tcircperiod} have their analogs here.

\section{Algorithm}

Let us consider any numerical $3$--semigroup $N=\sg{n_1,n_2,n_3}$ and $n\in N$. By Lemma~\ref{lem:coprim}, there is another semigroup $T=\sg{a,b,c}$, with $1\leq a<b<c$ and $\gcd(a,b)=\gcd(a,c)=\gcd(b,c)=1$, and $m\in T$ such that $\dd(n,N)=\dd(m,T)$. Lemma~\ref{lem:coprim} only requires a time cost of $O(\log n_3)$, in the worst case. Moreover, by Theorem~\ref{teo:EhII} and Theorem~\ref{teo:SeOz}, it can be assumed that $m\in\{0,\ldots,P-S\}$ with $P=abc$ and $S=a+b+c$.

There are three possible cases for the semigroup $T$,
\begin{itemize}
\item[(1)] $a>1$ and $c\notin\sg{a,b}$,
\item[(2)] $a>1$ and $c\in\sg{a,b}$,
\item[(3)] $a=1$.
\end{itemize}
Now we analyze each case for finding the related L-shapes. Then, the time cost of the related $S^\pm$ sums will be studied.

\subsection{Case 1: $a>1$ and $c\notin\sg{a,b}$}
In this case we have $e(T)=3$. From $c\notin\sg{a,b}$, we also have $c\leq\f(a,b)<(a-1)(b-1)<ab$.

\begin{lem}[{Rosales and Garc\'{\i}a-S\'anchez {\rm\cite[Chap.~9]{RoGa-NS:2009}}}]\label{lem:Mc}
Let $\sg{a,b,c}$ be a numerical $3$--se\-mi\-group with $1<a<b<c$. Assume that $\HH=\ele(l,h,w,y)$ is related to $T$ with $wy\neq0$. Then,
\[
M_c=(l-w)a+(h-y)b=(\delta+\theta)c=\min\{kc:~k\geq1, kc\in\sg{a,b}\}.
\]
\end{lem}
\noindent{\bf Proof}. Assume $k_0c<M_c$ for some $k_0\in\Nat$ with $k_0\geq1$. Then, $k_0c=\alpha a+\beta b$ with $\alpha,\beta\in\Nat$ and $\alpha+\beta\geq2$ (the identity $\alpha+\beta=1$ leads to $k_0c=a$ or $k_0c=b$, a contradiction to $a<b<c$).

Assume $\alpha\geq1$. Then, the squares $\sq{\alpha-1,\beta}$ and $\sq{l-w-1,h-y}$ represent the same equivalence class $[0]_c$. From $\sq{l-w-1,h-y}\in\HH$ (because of $y>0$) and $\sq{\alpha-1,\beta}\notin\HH$ (only one square in $\HH$ for each equivalence class), we have $(l-w-1)a+(h-y)b\leq(\alpha-1)a+\beta b$. Thus, $M_c\leq k_0c$ holds and makes a contradiction.

The case $\beta\geq1$ also makes a contradiction by similar arguments. $\square$

\begin{lem}\label{lem:cas1-unaL}
Let $T=\sg{a,b,c}$ be a numerical $3$-semigroup with $1<a<b<c$, $\gcd(a,b)=\gcd(a,c)=\gcd(b,c)=1$ and $c\notin\sg{a,b}$. Then, only one L-shape $\ele(l,h,w,y)$ is related to $T$ and $wy\neq0$.
\end{lem}
\noindent{\bf Proof}. Assume $w=0$. Then $hb=0+\theta c$ holds with $\theta\geq1$. Thus, $c\mid h$ ($\gcd(b,c)=1$) and $b\mid\theta$. Now, from $c=lh$, it follows that $h=c$ and $l=1$. So, $a=yb+\delta c$ holds and makes a contradiction. Indeed, either $\delta=0$ we have $y\neq0$ and $a\mid b$ or $\delta>0$ and $a\geq c$, a contradiction. The case $y=0$ also leads to contradiction by similar arguments.

According to \cite[theorems 2 and 3]{AM:2014}, $T$ has only one related L-shape $\HH=\ele(l,h,w,y)$ iff $(la-yb)(hb-wa)>0$. Assume $la=yb$ holds. Then, $a\mid y$ and $b\mid l$ ($\gcd(a,b)=1$) and $hb=wa+\theta c$ with $\theta\geq1$ (recall that $\delta+\theta\geq1$). From $c=lh-wy=\frac{yb}ah-wy=\frac ya(hb-wa)$, we have $ac=y\theta c$. So, $y\mid a$ holds and so $y=a$. Therefore, $\theta=1$ and $l=b$ hold.

Let us consider now $M_c$ defined in Lemma~\ref{lem:Mc}. Then,
\[
M_c=(l-w)a+(h-y)b=(b-w)a+(h-a)b=hb-wa=\theta c=c.
\]
So, $c\in\sg{a,b}$ holds and makes a contradiction. The assumption $hb=wa$ also leads to contradiction by similar arguments. $\square$

This lemma ensures that $yb<la$ and $wa<hb$ hold for $\ele(l,h,w,y)$ related to $T$. A direct consequence of Lemma~\ref{lem:cas1-unaL} is the non-symmetry of $T$.

\begin{lem}\label{lem:minimality-lh}
Let $T=\sg{a,b,c}$ be a numerical $3$-semigroup with $1<a<b<c$ and $\gcd(b,c)=\gcd(a,c)=1$. Assume $\HH=\ele(l,h,w,y)$ is an L-shape related to $T$. Then,
\begin{align*}
la=&\min\{ka:~k\geq1, ka\in\sg{b,c}\},\\
hb=&\min\{kb:~k\geq1, kb\in\sg{a,c}\}.
\end{align*}
\end{lem}
\noindent{\bf Proof}. Here we prove the first equality. The second one is proved by similar arguments.

As $\HH$ is an L-shape related to $T$, we have  $\Ap(c,T)=\{ia+jb:~\sq{i,j}\in\HH\}$. In particular, it follows that $la=\min\{ka:~ka\notin\Ap(c,T)\}$. Using the same notation of \cite[Lemma~10.18]{RoGa-NS:2009}, we have $c_1a=r_{12}b+r_{13}c$ with $r_{12},r_{13}>0$, where $c_1a=\min\{ka:~k\geq1,ka\in\sg{b,c}\}$.

As $la\notin\Ap(c,T)$, we have $la-c\in T$ and so $la-c=x_1a+x_2b+x_3c$. Assuming $x_1\neq0$, $(l-x_1)a-c=x_2b+x_3c\in T$ holds and then $(l-x_1)a\notin\Ap(c,T)$. This is a contradiction to the minimality of $la$. Therefore, $x_1=0$ holds. Thus, $la-c\in\sg{b,c}$ and $l\geq c_1$ from the minimality of $c_1a$.

Now, from $c_1a=r_{12}b+r_{13}c$ with $r_{13}>0$, it follows that $c_1a-c\in T$. Then, $c_1a\notin\Ap(c,T)$ and $l\leq c_1$ from the minimality of $la$. $\square$

\begin{lem}\label{lem:cas1-ub}
Let $T=\sg{a,b,c}$ be a numerical $3$-semigroup with $1<a<b<c$, $\gcd(a,b)=\gcd(a,c)=\gcd(b,c)=1$ and $c\notin\sg{a,b}$. Assume $T$ has only one related L-shape $\ele(l,h,w,y)$. Then, $h<a$ and $l<b$.
\end{lem}
\noindent{\bf Proof}. By Lemma~\ref{lem:minimality-lh}, it follows that $la\leq ab$ and $l\leq b$ holds. Similarly, $h\leq a$ also holds.

Assume $l=b$. So, $ab=la=yb+\delta c$ with $\delta\geq1$ (by Lemma~\ref{lem:cas1-unaL} we have $la>yb$). Then, $b(a-y)=\delta c$ holds and thus $c\mid(a-y)$ ($\gcd(b,c)=1$). That is, $a=y+\alpha c$ with $\alpha\geq1$ which contradicts inequality $a<c$. Similar arguments lead to contradiction assuming $h=a$. $\square$

In this case the sides of the L-shape are bounded by $w<l<b$ and $y<h<a$.

\subsection{Case 2: $a>1$ and $c\in\sg{a,b}$}

Identities $\gcd(a,b)=\gcd(b,c)=\gcd(a,c)=1$ also hold.

\begin{lem}[A. and Mariju\'an {\rm \cite[Theorem~8-(d)]{AM:2014}}]\label{lem:cas2}
Assume $c=\lambda a+\mu b$, $\lambda,\mu\in\Nat$, $0<\mu<a$, $\gcd(a,\mu)=\gcd(b,\lambda)=1$. Then, $\lambda\neq b$ and there are two L-shapes related to $T=\sg{a,b,c}$, $\HH_1=\ele(\lambda+b,a,b,a-\mu)$ with $(\delta,\theta)=(1,0)$ and $\HH_2$ with $(\delta,\theta)=(0,1)$ given by
\[
\HH_2=\begin{cases}
\ele(b,a+\mu,b-\lambda,a)&\text{if }\lambda<b,\\
\ele(b,(1+\lfloor\lambda/b\rfloor)a+\mu,b-s,a)&\text{if }\lambda>b\text{ where }\lambda=\lfloor\lambda/b\rfloor b+s, 0\leq s<b.
\end{cases}
\]
\end{lem}

\subsection{Case 3: $a=1$}

Consider a semigroup $T=\sg{1,b,c}$ with $1<b<c$ and $\gcd(b,c)=1$.

\begin{lem}\label{lem:cas3}
Consider the numerical semigroup $T=\sg{1,b,c}$ with $\gcd(b,c)=1$. Then, there are two related L-shapes $\HH_1=\ele(c,1,b,0)$ with parameters $(\delta,\theta)=(1,0)$ and $\HH_2$ with parameters $(\delta,\theta)=(0,1)$ given by
\[
\HH_2=\begin{cases}
\ele(b,2,2b-c,1)&\text{if }c<2b,\\
\ele(b,1+\lfloor c/b\rfloor,b-r,1)&\text{if }c>2b\text{ where }c=\lfloor c/b\rfloor b+r, 0\leq r<b.
\end{cases}
\]
\end{lem}
\noindent{\bf Proof}. $\HH_1$ is related to $T$ by Theorem~\ref{teo:LMDD}. As $\theta=0$, using the transformation of L-shapes defined in \cite[Theorem~3]{AM:2014}, we obtain $\HH_2$ from $\HH_1$. $\square$

\section{Time cost}

Let us analyze now the time cost, in the worst case. This analysis will be done under the assumption of $m\approx P=abc$. This is the same assumption as the one made in the analysis of algorithms \algP, \algL\ and \algBCS.

Applying Theorem~\ref{teo:casi}, Theorem~\ref{teo:casii} or Theorem~\ref{teo:casiii} requires the calculation of the L-shape $\HH=\ele(l,h,w,y)$, the related basic factorization $(x_0,y_0,z_0)$ and all the related $S^\pm$ sums. The first two calculations have a time cost of $O(\log c)$ \cite{AB:2008}. Then, all $S^\pm$ have to be calculated.

Consider a generic sum $S^\pm(s,t,q,N)=\sum_{k=0}^N\left\lfloor\frac{s\pm kt}q\right\rfloor$, with $0\leq s,t<q$. Using the same notation of Section~\ref{sec:Spm}, we have
\begin{itemize}
\item If $t\mid q$, Theorem~\ref{teo:sm-1} for $S^+$ and Theorem~\ref{teo:smenys-1} for $S^-$ ensure a constant time cost.
\item If $t\nmid q$ and $\gcd(t,q)=1$, Theorem~\ref{teo:sm-2} for $S^+$ or Theorem~\ref{teo:smenys-2} for $S^-$ has to be applied. All computations are focused on finding the subsets of hS indices $K_1\subset\{0,\ldots,M-1\}$, $K_2\subset\{j0+ut,\ldots, M-1\}$ or $J\subset\{0,\ldots,t-1\}$. As $|K_1|,|K_2|\leq|J|$, the cost has order $O(\widehat{q})\leq O(t)$, in the worst case.
\item If $t\nmid q$ and $\gcd(t,q)=g>1$, consider $\tilde{t}=t/g$. Then, Theorem~\ref{teo:sm-3} or Theorem~\ref{teo:smenys-3} and similar arguments as in the previous case ensure a time cost upper bounded by $O(\tilde{t})$.
\end{itemize}

\begin{rem}\sl\label{rem:ordre}
Previous comments point to the fact that the higher cost of computation is reached when $t\nmid q$ and $\gcd(t,q)=1$. In this case, the time cost is upperbounded by $O(t)$.
\end{rem}

Given a semigroup $S=\sg{n_1,n_2,n_3}$ and $n\in S$, apply Lemma~\ref{lem:coprim} at constant time cost for obtaining $T=\sg{a,b,c}$ with $1\leq a<b<c$ and $\gcd(a,b)=\gcd(b,c)=\gcd(a,c)=1$ and $m\in T$ such that $\dd(n,S)$ can be calculated from $\dd(m,T)$. Then, we have to analyze the time cost of each case given in the previous section.

The worst case for the calculation of $S^\pm(s,t,q,N)$, as it is highlighted in Remark~\ref{rem:ordre}, appears when $t\nmid q$ and $\gcd(t,q)=1$. This case will be assumed in all cases in the following analysis. Thus, the resulting worst case order will be a pessimistic estimation.

\begin{itemize}
\item Case 1 ($a>1$, $c\notin\sg{a,b}$). By Lemma~\ref{lem:cas1-unaL}, the L-shape $\HH$ belongs to the case (iii) $\delta\theta>0$, subcase (iii.4) $wy\neq0$. The following sums have to be evaluated
\begin{itemize}
\item If $k_1=0$, there is one sum $S^-_1=\sum_{k=0}^{A_m}\left\lfloor\frac{\widehat{z_{0,1}}-k\widehat{\theta}}{\delta}\right\rfloor$. As $0\leq\widehat{\theta}<\delta=\frac{la-yb}c<\frac{la}c<\frac{ab}c<a$, the cost of calculating $S_1^-$ is upperbounded by $O(a)$.
\item If $1\leq k_1\leq A_m$, there are two sums $S^+_1=\sum_{k=0}^{k_1-1}\left\lfloor\frac{\widehat{y_0}+k\widehat{h}}{y}\right\rfloor$ and $S^-_1=\sum_{k=0}^{A_m}\left\lfloor\frac{\widehat{z_{0,1}}-k\widehat{\theta}}{\delta}\right\rfloor$. The calculation of $S^-_1$ is $O(a)$. As $\widehat{h}<y<a$, the order for calculating $S^+_1$ is also upperbounded by $O(a)$. Thus, the worst case order of this case is $O(a)$.
\item If $k_1>A_m$, there is one sum $S_2^+=\sum_{k=0}^{A_m}\left\lfloor\frac{\widehat{y_0}+k\widehat{h}}{y}\right\rfloor$. This sum has the same order as $S_1^+$, that is $O(a)$.
\item If $k_0=0$, there is one sum $S_2^-=\sum_{k=0}^{A_m}\left\lfloor\frac{\widehat{z_{0,2}}-k\widehat{\delta}}{\theta}\right\rfloor$. From $\widehat{\delta}<\theta=\frac{hb-wa}c<\frac{hb}c<\frac{ab}c<a$, the order is upperbounded by $O(a)$.
\item If $1\leq k_0\leq A_m$, there are two sums $S_3^+=\sum_{k=0}^{k_0-1}\left\lfloor\frac{\widehat{x_0}+k\widehat{l}}{w}\right\rfloor$ and $S^-_2$. From $\widehat{l}<w<b$, the cost of both calculations is $O(a)+O(b)=O(b)$.
\item If $k_0>A_m$, we have $S_4^+=\sum_{k=0}^{A_m}\left\lfloor\frac{\widehat{x_0}+k\widehat{l}}{w}\right\rfloor$ with the same cost of $S_3^+$, that is $O(b)$.
\end{itemize}
Therefore, the overall cost of the Case~1 is $O(b)$.

\item Case 2 ($a>1$, $c\in\sg{a,b}$). Let us consider $c=\lambda a+\mu b$ with $1\leq\mu<a$. By Lemma~\ref{lem:cas2} there are three possible cases to be examined.

Consider the L-shape $\HH_1=\ele(\lambda+b,a,b,a-\mu)$ with $\delta=1$ and $\theta=0$. Then, $A_m=z_0$ and $k_1=\left\lceil\frac{z_0a-y_0}a\right\rceil\leq z_0=A_m$. Thus, the case $k_1>A_m$ never appears. So,
\begin{itemize}
\item If $k_1=0$, there is one sum $S_1^+=\sum_{k=0}^{A_m}\left\lfloor\frac{\widehat{x_0}+k\widehat{l}}{w}\right\rfloor$ with $w=b$. Then, the order is upperbounded by $O(b)$.
\item If $1\leq k_1\leq A_m$, there are two sums $S_1^+$ and $S_2^+=\sum_{k=0}^{k_1-1}\left\lfloor\frac{\widehat{y_0}+k\widehat{h}}{y}\right\rfloor$, From $y=a-\mu<a$, the order is upperbounded by $O(a)+O(b)=O(b)$.
\end{itemize}
Let us examine the other related L-shape $\HH_2$ which have an expression depending on $\lambda$. Assume $\lambda<b$\label{abcd}. Then, we have $\HH_2=\ele(b,a+\mu,b-\lambda,a)$ with $\delta=0$ and $\theta=1$. As $\delta=0$, this is the Case-(i) with $w=b-\lambda>0$. So, $A_m=z_0$ holds and the case $k_0>A_m$ never appears. Then,
\begin{itemize}
\item If $k_0=0$, there is one sum $S_1^+=\sum_{k=0}^{A_m}\left\lfloor\frac{\widehat{y_0}+k\widehat{h}}{y}\right\rfloor$. From $y=a$, we have order upperbounded by $O(a)$.
\item If $1\leq k_0\leq A_m$, there are two sums $S_1^+$ and $S_2^+=\sum_{k=0}^{k_0-1}\left\lfloor\frac{\widehat{x_0}+k\widehat{l}}{w}\right\rfloor$. From $w=b-\lambda<b$, the order is upperbounded by $O(a)+O(b)=O(b)$.
\end{itemize}
Assume now $\lambda>b$. Then, the related L-shape is $\HH_2=\ele(b,(1+\lfloor\lambda/b\rfloor)a+\mu,b-s,a)$ with $0\leq s<b$, $\delta=0$ and $\theta=1$. We also are in the Case (i) with $w=b-s>0$. So, $A_m=z_0$ and $k_0\leq A_m$ always holds. Then,
\begin{itemize}
\item If $k_0=0$, there is one sum $S_1^+=\sum_{k=0}^{A_m}\left\lfloor\frac{\widehat{y_0}+k\widehat{h}}{y}\right\rfloor$. From $y=a$, the order is upperbounded by $O(a)$.
\item If $1\leq k_0\leq A_m$, there are two sums $S_1^+$ and $S_2^+=\sum_{k=0}^{k_0-1}\left\lfloor\frac{\widehat{x_0}+k\widehat{l}}{w}\right\rfloor$. The order is upperbounded by $O(a)+O(b)=O(b)$.
\end{itemize}
In any case, using either $\HH_1$ or $\HH_2$, the overall order is upperbounded by $O(b)$.

\item Case 3 ($a=1$). We have $\gcd(b,c)=1$. By Lemma~\ref{lem:cas3}, there are three possibilities to analyze.

Consider $\HH_1=\ele(c,1,b,0)$, with $\delta=1$ and $\theta=0$. This L-shape can be used in the two cases $c<2b$ and $c>2b$ (note that $c\neq2b$). Look at the Case (ii) in Section~\ref{sec:casii} and Theorem~\ref{teo:casii}. As $y=0$, from \eqref{eq:casii1} we have to calculate one sum $S_1^+=\sum_{k=0}^{A_m}\left\lfloor\frac{\widehat{x_0}+k\widehat{l}}{w}\right\rfloor$. From $\widehat{l}<w=b$, the order is upperbounded by $O(b)$.

Let us consider now the case $c<2b$. Again by Lemma~\ref{lem:cas3}, we can use the L-shape $\HH_2=\ele(b,2,2b-c,1)$ with $\delta=0$ and $\theta=1$. Note that $A_m=z_0$. We have to look at Theorem~\ref{teo:casi}-(i.2) ($y=1>0$). We have $k_0\leq A_m$, then the case (i.2.3) never appears. Then,
\begin{itemize}
\item If $k_0=0$, there is only one sum given by \eqref{eq:casi1} $S_1^+=\sum_{k=0}^{A_m}\left\lfloor\frac{\widehat{y_0}+k\widehat{h}}{y}\right\rfloor$. From $y=1$, the cost is constant.
\item If $1\leq k_0\leq A_m$, there are two sums $S_1^+$ and $S_2^+=\sum_{k=0}^{k_0-1}\left\lfloor\frac{\widehat{x_0}+k\widehat{l}}{w}\right\rfloor$. From $\widehat{l}<w=2b-c<b$, the cost is upperbounded by $O(b)$. Then, the total cost is upperbounded by $O(b)+O(1)=O(b)$.
\end{itemize}
When $c>2b$, we can use the L-shape $\HH_2=\ele(b,1+\lfloor c/b\rfloor,b-r,1)$ with $c=\lfloor c/b\rfloor b+r$ and $0\leq r<b$. The related parameters are $\delta=0$ and $\theta=1$. Using the same arguments of the previous case, it follows that the overall order is upperbounded by $O(b)$.

Therefore, using any admissible L-shape, the total cost of this case is also $O(b)$.
\end{itemize}

So, we need a cost of $O(\log c)$ to calculate the related L-shape and the basic coordinates plus $O(b)$ to calculate the involved $S^\pm$ sums. Hence, our algorithm has a time cost of $O(b+\log c)$. Common instances of semigroups are such that $O(\log c)\ll O(b)$. Then, we have the following result.

\begin{teo}\label{teo:cost}
The time cost, in the worst case, for computing the denumerant $\dd(m,T)$ is upperbounded by $O(b+\log c)$.
\end{teo}

\begin{rem}\sl
When many instances of $m\in T$ are given and the semigrup $T$ is fixed, the related L-shape is computed only once. Thus, the first calculated denumerant has a time cost of $O(b+\log c)$. The subsequent instances only need a time cost of $O(b)$.
\end{rem}

\section{Some time tests}

All the computations of this section have been made using SageMath~7.3 \cite{Sage} and non compiled code on a \textsf{i5@1.3Ghz} processor. Here we test our algorithm, denoted by \algAL, versus the algorithms \algP, \algL\ and \algBCS. In the following, we use the notation $P=abc$ and $S=a+b+c$. The time required to calculate denumerants highly depends on the selected semigroup. This fact is reflected in the following subsections. All semigroups in this section will meet the property $\gcd(a,b)=\gcd(a,c)=\gcd(b,c)=1$. By Lemma~\ref{lem:coprim} of Brown, Chou and Shiue, this restrictions does not represent any loose of generality.

Time costs of the involved algorithms are by Table~\ref{tab:tc}. It is assumed that $T=\sg{a,b,c}$ and $m\approx P=abc$.
\begin{table}[h]
\centering
%\footnotesize
%\small
\begin{tabular}{|c|l|}
\hline
Algorithm&Time cost\\
\hline\hline
\algP&$O(c\log c)$\\
\algL&$O(ab\log b)$\\
\algBCS&$O(ab\log c)$\\
\algAL&$O(b+\log c)$\\
\hline
\end{tabular}
\caption{Time costs for $T=\sg{a,b,c}$ and $m\approx P$}
\label{tab:tc}
\end{table}

\begin{rem}\sl\label{rem:tc}
According to Table~\ref{tab:tc} there are some generic behaviours to be highlighted:
\begin{itemize}
\item[(i)] Algorithm~\algAL\ has the best time cost.
\item[(ii)] When $a=1$, algorithms \algL\ and \algBCS\ are faster than Algorithm~\algP\ when $b\ll c$. However, when $c\approx b$, algorithms \algP, \algL\ and \algBCS\ run at similar speed.
\item[(iii)] When $a\neq1$, there are two different behaviours,
\begin{itemize}
\item[(iii.1)] if $ab<c$, Algorithm~\algL\ is faster than algorithms \algBCS\ and \algP,
\item[(iii.2)] otherwise, when $ab>c$, Algorithm~\algP\ wins \algL\ and \algBCS. 
\end{itemize}
\item[(iv)] When $b\approx c$, Algorithm~\algP\ is faster than algorithms \algL\ and \algBCS\ provided that $a\gg1$.
\end{itemize}
\end{rem}

In the following subsections we take elements $m$ of the semigroup that are closed to $P-S$.

\subsection{$a>1$, $c\notin\sg{a,b}$}

\begin{table}[h]
\centering
%\footnotesize
\small
\begin{tabular}{|r|r|r|r|r|r|r|}
\hline
$k$&\multicolumn{1}{|c}{$m_k$}&\multicolumn{1}{|c}{$\dd(m_k,T_{1,k})$}&\multicolumn{1}{|c}{\algP} &\multicolumn{1}{|c}{\algL} &\multicolumn{1}{|c}{\algBCS}&\multicolumn{1}{|c|}{\algAL}\\
\hline\hline
1&4465&2232&0.002304&0.003423&0.010860&0.000291\\
2&34139180&17069589&0.078596&0.116350&1.025509&0.000615\\
3&207657687311&103828843654&5.291058&9.787089&68.891171&0.000533\\
4&1235137178269914&617568589134955&424.791713&740.592501&5275.727091&0.000376\\
\hline
\end{tabular}
\caption{$T_{1,k}=\sg{7^k,11^k,\f(7^k,11^k)}$, $m_k=P_k-S_k-k$}
\label{tab:7-11-fn}
\end{table}

\begin{table}[h]
\centering
%\footnotesize
\small
\begin{tabular}{|r|r|r|r|r|r|r|}
\hline
$k$&\multicolumn{1}{|c}{$m_k$}&\multicolumn{1}{|c}{$\dd(m_k,T_{2,k})$}&\multicolumn{1}{|c}{\algP} &\multicolumn{1}{|c}{\algL} &\multicolumn{1}{|c}{\algBCS}&\multicolumn{1}{|c|}{\algAL}\\
\hline\hline
1&893&446&0.001188&0.003566&0.006692&0.000545\\
2&723044&361521&0.005573&0.152317&0.419318&0.000326\\
3&608098947&304049472&0.045436&10.102134&31.251389&0.000402\\
\hline
\end{tabular}
\caption{$T_{2,k}=\sg{7^k,11^k,11^k+1}$, $m_k=P_k-S_k-k$}
\label{tab:7-11-11p1}
\end{table}

In this case, inequality $c<ab$ always holds. Then, as it has been comment in Remark~\ref{rem:tc}-(iii.2), Algorithm~\algP\ is  faster than Algorithm~\algL\ and Algorithm~\algBCS. Table~\ref{tab:7-11-fn} shows how this assertion is kept for the semigroups $T_{1,k}=\sg{7^k,11^k,\f(7^k,11^k)}$ and $m_k=P_k-S_k-k$ for $k\in\{1,2,3,4\}$. The non increasing sequence of times in the column of Algorithm~\algAL\ is because the corresponding L-shapes. Their entries do not always increase as the value of $k$ does.

Table~\ref{tab:7-11-11p1}, for the semigroups $T_{2,k}=\sg{7^k,11^k,11^k+1}$, shows an instance of the case $c\approx b$ and, as it has been noticed in Remark~\ref{rem:tc}-(iv), Algorithm~\algP\ is faster than algorithms \algL\ and \algBCS.

\subsection{$a>1$, $c\in\sg{a,b}$}

In this subsection we take the semigroups $T_{3,k}=\sg{7^k,11^k,7^{k}+11^{2k}}$ for the case $ab<c$ and $T_{4,k}=\sg{7^k,11^k,7^k+11^k}$ for $ab>c$. Tables \ref{tab:7-11-a2bp1} and \ref{tab:7-11-apb} show the influence of inequalities $ab<c$ and $ab>c$ in the resulting time cost. Here, item (iii) of Remark~\ref{rem:tc} is also clear.

\begin{table}[h]
\centering
%\footnotesize
\small
\begin{tabular}{|r|r|r|r|r|r|r|}
\hline
$k$&\multicolumn{1}{|c}{$m_k$}&\multicolumn{1}{|c}{$\dd(m_k,T_{3,k})$}&\multicolumn{1}{|c}{\algP} &\multicolumn{1}{|c}{\algL} &\multicolumn{1}{|c}{\algBCS}&\multicolumn{1}{|c|}{\algAL}\\
\hline\hline
1&9709&4854&0.004208&0.002923&0.015249&0.000422\\
2&87082148&43541073&0.188668&0.133961&1.012403&0.000332\\
3&808930875251&404465437624&19.829875&9.667596&69.573355&0.000392\\
\hline
\end{tabular}
\caption{$T_{3,k}=\sg{7^k,11^k,7^{k}+11^{2k}}$, $m_k=P_k-S_k-k$}
\label{tab:7-11-a2bp1}
\end{table}

\begin{table}[h]
\centering
%\footnotesize
\small
\begin{tabular}{|r|r|r|r|r|r|r|}
\hline
$k$&\multicolumn{1}{|c}{$m_k$}&\multicolumn{1}{|c}{$\dd(m_k,T_{4,k})$}&\multicolumn{1}{|c}{\algP} &\multicolumn{1}{|c}{\algL} &\multicolumn{1}{|c}{\algBCS}&\multicolumn{1}{|c|}{\algAL}\\
\hline\hline
1&1349&674&0.001330&0.002596&0.009378&0.000342\\
2&1007588&503793&0.006507&0.155632&0.489768&0.000453\\
3&764232891&382116444&0.045628&9.864927&35.554438&0.000263\\
\hline
\end{tabular}
\caption{$T_{4,k}=\sg{7^k,11^k,7^k+11^k}$, $m_k=P_k-S_k-k$}
\label{tab:7-11-apb}
\end{table}

\subsection{$a=1$}
Let us take the semigroups $T_{5,k}=\sg{1,7^k,11^k}$. Table~\ref{tab:1-7-11} confirms that Algorithm~\algP\ is slower than algorithms \algL\ and \algBCS. This rule is not noticeable with respect to Algorithm~\algBCS\ for small values of $k$. However, it turns apparent from the value $k=6$.

\begin{table}[h]
\centering
%\footnotesize
\small
\begin{tabular}{|r|r|r|r|r|r|r|}
\hline
$k$&$m_k$&$\dd(m_k,T_{5,k})$&\algP &\algL &\algBCS&\algAL\\
\hline\hline
1&57&29&0.001047&0.000917&0.002000&0.000381\\
2&5756&2878&0.003230&0.002486&0.006424&0.000301\\
3&454855&227427&0.024524&0.009076&0.056042&0.000251\\
4&35135994&17567996&0.212520&0.048399&0.383123&0.000379\\
5&2706606293&1353303145&2.070613&0.350236&2.447490&0.001594\\
6&208420490872&104210245434&20.946929&2.581581&16.702903&0.001885\\
7&16048502956131&8024251478063&225.232675&16.902235&116.667887&0.009188\\
\hline
\end{tabular}
\caption{$T_{5,k}=\sg{1,7^k,11^k}$, $m_k=P_k-S_k-k$}
\label{tab:1-7-11}
\end{table}

\begin{table}[h]
\centering
%\footnotesize
\small
\begin{tabular}{|r|r|r|r|r|r|r|}
\hline
$k$&$m_k$&$\dd(m_k,T_{6,k})$&\algP &\algL &\algBCS&\algAL\\
\hline\hline
1&39&20&0.000963&0.000689&0.001357&0.000307\\
2&2348&1174&0.002826&0.001878&0.006299&0.000307\\
3&117301&58650&0.013340&0.008977&0.033656&0.000405\\
4&5762394&2881196&0.065825&0.061175&0.203000&0.000253\\
5&282458435&141229216&0.382746&0.397065&1.322981&0.000200\\
6&13841169544&6920584770&2.668743&2.589575&9.066692&0.000212\\
7&678222249297&339111124646&17.936306&17.299389&61.854215&0.000212\\
8&33232924804790&16616462402392&127.708571&120.464809&439.603427&0.000302\\
\hline
\end{tabular}
\caption{$T_{6,k}=\sg{1,7^k,7^k+1}$, $m_k=P_k-S_k-k$}
\label{tab:a1bis}
\end{table}

Consider now the semigroups $T_{6,k}=\sg{1,7^k,7^k+1}$, where $c\approx b$. According to Remark~\ref{rem:tc}-(ii), the algorithms \algP, \algL\ and \algBCS\ have similar time cost. Table~\ref{tab:a1bis} shows this behaviour in algorithms \algP\ and \algL. Algorithm~\algBCS\ runs between three and four times slower.

\subsection{Almost medium and large input data}
Now we take larger input values for the Algorithm~\algAL. Usually, our algorithm can manage almost middle input values at acceptable time output. However, when the involved $S^\pm$ sums take some proper parameters, the time cost can be almost constant. These cases allow the Algorithm~\algAL\ to take large input values.

We consider the same semigroups of the previous sections to see these behaviours. When the output values $m_k$ and $\dd(m_k,T)$ turn to be large, tables will show $\ell(m_k)$ and $\ell(\dd(m_k,T))$.

Table~\ref{tab:T1-large} and Table~\ref{tab:T2-large} belong to the case $a>1$ with $c\notin\sg{a,b}$. The case $a>1$ with $c\in\sg{a,b}$, is represented by Table~\ref{tab:T3-large} when $ab<c$ and Table\ref{tab:T4-large} when $ab>c$. Finally, the case $a=1$ is represented by tables \ref{tab:T5-large} and \ref{tab:T6-large}.

\begin{table}[h]
\centering
%\footnotesize
\small
\begin{tabular}{|r|r|r|r|}
\hline
$k$&$\ell(m_k)$&$\ell(\dd(m_k,T_{1,k}))$& \multicolumn{1}{c|}{\algAL}\\
\hline\hline
10&38&38&0.000424\\
100&378&377&0.000966\\
1000&3773&3773&0.002669\\
10000&37730&37730&0.051687\\
100000&377299&377298&1.350448\\
1000000&3772982&3772982&25.325951\\
\hline
\end{tabular}
\caption{$T_{1,k}=\sg{7^k,11^k,\f(7^k,11^k)}$, $m_k=P_k-S_k-k$}
\label{tab:T1-large}
\end{table}

\begin{table}[h]
%\centering
\hspace*{-7mm}
\footnotesize
%\small
\begin{tabular}{|r|r|r|r|}
\hline
$k$&\multicolumn{1}{|c}{$m_k$}&\multicolumn{1}{|c}{$\dd(m_k,T_{2,k})$}& \multicolumn{1}{|c|}{\algAL}\\
\hline\hline
4&514710794634&257355397315&0.000495\\
5&435933001714249&217966500857122&0.000889\\
6&369233168511568240&184616584255784117&0.000775\\
7&312740333247126511823&156370166623563255908&0.003953\\
8&264891049902986514370070&132445524951493257185031&0.010668\\
9&224362718316312996430224405&112181359158156498215112198&0.007589\\
10&190035222340650307226923642236&95017611170325153613461821113&0.013501\\
11&160959833316889266300917603625499&80479916658444633150458801812744&0.155148\\
12&136332978818970809672136276502054306&68166489409485404836068138251027147&5.684347\\
13&115474033059634827078142773316758235361&57737016529817413539071386658379117674&0.888997\\
14&97806506001508122984196519581781213521032&48903253000754061492098259790890606760509&31.661762\\
15&82842110583277181850188186745272781621519975&41421055291638590925094093372636390810759980&80.498716\\
\hline
\end{tabular}
\caption{$T_{2,k}=\sg{7^k,11^k,11^k+1}$, $m_k=P_k-S_k-k$}
\label{tab:T2-large}
\end{table}

Algorithm~\algAL\ allows almost middle length inputs, above a hundred digits. Several instances of this inputs at reasonable time output are given in tables \ref{tab:T2-large}, \ref{tab:T3-large} and \ref{tab:T5-large}. The nature of the involved $S^\pm$ sums has an interesting property. Some parameters taken by these sums make almost constant the time cost of the denumerant's calculation. In these cases, the algorithm can handle large inputs (million digits) at a small time cost. Tables \ref{tab:T1-large}, \ref{tab:T4-large} and \ref{tab:T6-large} show some instances of this good behaviour.

\begin{table}[h]
%\centering
\hspace*{-9mm}
\footnotesize
%\small
\begin{tabular}{|r|r|r|r|}
\hline
$k$&\multicolumn{1}{|c}{$m_k$}&\multicolumn{1}{|c}{$\dd(m_k,T_{3,k})$}& \multicolumn{1}{|c|}{\algAL}\\
\hline\hline
4&7535450720580234&3767725360290115&0.000684\\
5&70207055450352553785&35103527725176276890&0.009376\\
6&654118736532593706215344&327059368266296853107669&0.003036\\
7&6094424053060467191130813247&3047212026530233595565406620&0.308053\\
8&56781748786352120105926189224470&28390874393176060052963094612231&0.497334\\
9&529035553379910639083032546370845061&264517776689955319541516273185422526&5.588338\\
10&4929024250806922465243407641240300023740&2464512125403461232621703820620150011865&52.150354\\
11&45923718944749929614349668788645626831480843&22961859472374964807174834394322813415740416&248.336705\\
\hline
\end{tabular}
\caption{$T_{3,k}=\sg{7^k,11^k,7^k+11^{2k}}$, $m_k=P_k-S_k-k$}
\label{tab:T3-large}
\end{table}

\begin{table}[h]
\centering
%\footnotesize
\small
\begin{tabular}{|r|r|r|r|}
\hline
$k$&$\ell(m_k)$&$\ell(\dd(m_k,T_{4,k}))$& \multicolumn{1}{c|}{\algAL}\\
\hline\hline
10&30&29&0.000494\\
100&293&293&0.000435\\
1000&2928&2928&0.002186\\
10000&29279&29279&0.046576\\
100000&292789&292789&1.200240\\
1000000&2927884&2927884&23.496686\\
\hline
\end{tabular}
\caption{$T_{4,k}=\sg{7^k,11^k,7^k+11^k}$, $m_k=P_k-S_k-k$}
\label{tab:T4-large}
\end{table}

Now, we briefly comment this almost constant time cost behaviour of Algorithm~\algAL\ in tables \ref{tab:T1-large}, \ref{tab:T4-large} and \ref{tab:T6-large}. In fact, almost all the time is spent in the computation of the related L-shape, that is $O(\log c)$.

The semigroup $T_{1,n}=\sg{7^n,11^n,\f(7^n,11^n)}$, from Theorem~\ref{teo:LMDD}, has related the L-shape $\HH_{1,n}=\ele(11^n-1,7^n-1,1,1)$ with $\delta=\theta=1$. Then, this is the case $a>1$ with $c\notin\sg{a,b}$. We have to calculate some of the sums $\sum_{k=0}^{A_m}\left\lfloor\frac{\widehat{z_{0,1}}-k\widehat{\theta}}{\delta}\right\rfloor$, $\sum_{k=0}^{A_m}\left\lfloor\frac{\widehat{z_{0,2}}-k\widehat{\delta}}{\theta}\right\rfloor$, $\sum_{k=0}^{k_1-1}\left\lfloor\frac{\widehat{y_0}+k\widehat{h}}{y}\right\rfloor$, $\sum_{k=0}^{A_m}\left\lfloor\frac{\widehat{y_0}+k\widehat{h}}{y}\right\rfloor$, $\sum_{k=0}^{k_0-1}\left\lfloor\frac{\widehat{x_0}+k\widehat{l}}{w}\right\rfloor$ and  $\sum_{k=0}^{k_0-1}\left\lfloor\frac{\widehat{x_0}+k\widehat{l}}{w}\right\rfloor$. From $\delta=\theta=w=y=1$, all these sums are calculated at constant time. Therefore, the fast computation of denumerants in Table~\ref{tab:T1-large} is clear now.

\begin{table}[h]
\centering
%\hspace*{-8mm}
%\footnotesize
\small
\begin{tabular}{|r|r|r|r|}
\hline
$k$&\multicolumn{1}{|c}{$m_k$}&\multicolumn{1}{|c}{$\dd(m_k,T_{5,k})$}& \multicolumn{1}{|c|}{\algAL}\\
\hline\hline
8&1235736071423990&617868035711992&0.465418\\
9&95151692050870129&47575846025435061&5.482949\\
10&7326680446366300788&3663340223183150390&18.740051\\
11&564154396101848452607&282077198050924226299&241.083486\\
\hline
\end{tabular}
\caption{$T_{5,k}=\sg{1,7^k,11^{k}}$, $m_k=P_k-S_k-k$}
\label{tab:T5-large}
\end{table}

\begin{table}[h]
\centering
%\footnotesize
\small
\begin{tabular}{|r|r|r|r|}
\hline
$k$&$\ell(m_k)$&$\ell(\dd(m_k,T_{6,k}))$& \multicolumn{1}{c|}{\algAL}\\
\hline\hline
10&17&17&0.000623\\
100&170&169&0.000652\\
1000&1691&1690&0.000622\\
10000&16902&16902&0.003141\\
100000&169020&169020&0.061677\\
1000000&1690197&1690196&0.742673\\
10000000&16901961&16901961&10.886006\\
\hline
\end{tabular}
\caption{$T_{6,k}=\sg{1,7^k,7^k+1}$, $m_k=P_k-S_k-k$}
\label{tab:T6-large}
\end{table}

The semigroup $T_{4,n}=\sg{7^n,11^n,7^n+11^n}$ has related the L-shape $\HH_{4,n}=\ele(11^n,7^n+1,11^n-1,7^n)$ with $\delta=0$ and $\theta=1$. This is the case $a>1$ with $c\in\sg{a,b}$ and parameters $\lambda=\mu=1$ and $\lambda<b=11^n$. Thus, following this case at page~\pageref{abcd}, we have two possibilities:
\begin{itemize}
\item When $k_0=0$, it has to be computed the sum $\sum_{k=0}^{A_m}\left\lfloor\frac{\widehat{y_0}+k\widehat{h}}{y}\right\rfloor$ with $h=7^n+1$ and $y=7^n$. From $h=y+1$, it follows that $\widehat{h}=1$ and the sum can be computed at constant cost from Theorem~\ref{teo:sm-1}.
\item Otherwise, when $1\leq k_0\leq A_m$, the algorithm calculates the sum $\sum_{k=0}^{k_0-1}\left\lfloor\frac{\widehat{x_0}+k\widehat{l}}{w}\right\rfloor$. Then, $l=w+1$ holds and, by the previous argument, the sum can be computed at constant time cost.
\end{itemize}
Therefore, the fast behaviour of the algorithm in Table~\ref{tab:T4-large} is now clear.

Finally, let us consider the semigroups $T_{6,n}=\sg{1,7^n,7^n+1}$. A related L-shape is $\HH_{6,n}=\ele(7^n,2,7^n-1,1)$ with $\delta=0$ and $\theta=1$. This is the case $a=1$ with parameters $\lambda=\mu=1$ and $c<2b$. Here we also have two possible cases:
\begin{itemize}
\item When $k_0=0$, there is only one sum to be computed, $S_1^+=\sum_{k=0}^{A_m}\left\lfloor\frac{\widehat{y_0}+k\widehat{h}}{y}\right\rfloor$. Here we have  $y=1$. Thus, this sum is calculated at constant time.
\item If $1\leq k_0\leq A_m$, we have to compute $S_1^+$ of the previous case and $S_2^+=\sum_{k=0}^{k_0-1}\left\lfloor\frac{\widehat{x_0}+k\widehat{l}}{w}\right\rfloor$. Again, from $l=w+1$ and Theorem~\ref{teo:sm-1}, the sum $S_2^+$ can also be calculated at constant time cost.
\end{itemize}
Thus, the speed of the algorithm in Table~\ref{tab:T6-large} is now clear.

\begin{rem}\label{rem:fast}\sl
Many semigroups have related an L-shape $\ele(l,h,w,y)$ with $\delta=1$ and/or $\theta=1$, $w=1$ and/or $y=1$. Additionally, many elements of the semigroup $m\in T$ have null coefficient multiplying $k$ in the $S^\pm$ sums. So, the fast behaviour of this algorithm eventually can be habitual.
\end{rem}

\section{Conclusion}
Algorithm~\algAL\ accepts almost medium input data to calculate denumerants of numerical $3$-semigroups at acceptable speed using an ordinary computer (tables \ref{tab:T2-large}, \ref{tab:T3-large} and \ref{tab:T5-large}). As far as we know, this algorithm is faster than usual known implemented algorithms for embedding dimension three numerical semigroups. This is the behaviour in the worst case. Eventually, this algorithm accepts large input data (tables \ref{tab:T1-large}, \ref{tab:T4-large} and \ref{tab:T6-large}).

The main tool of this algorithm is the hS-type set of ordered indices of intervals. As the computation techniques for obtaining these sets become faster, the time cost of this algorithm turns to be smaller.

It is difficult to generalize the algorithm to larger embedding dimensions because of the related minimum distance diagrams. Less is known about these diagrams related to numerical $n$-semigroups for $n\geq4$, mainly a generic geometrical description.

\end{document}